\documentclass[10pt,twocolumn,twoside]{IEEEtran}

\usepackage{ifpdf, flushend,subfigure}
\usepackage{graphicx}
\usepackage{epstopdf}
\graphicspath{{../}}
\usepackage{multirow}
\usepackage{float}
\usepackage[cmex10]{amsmath}
\usepackage {amssymb}
\usepackage{array}
\usepackage{mdwmath}
\usepackage{mdwtab}
\usepackage{eqparbox}
\usepackage{nomencl}
\usepackage{cite}
\usepackage{algorithm}
\usepackage{algorithmic}
\usepackage{makeidx}
\usepackage{ifthen}
\usepackage{subfigure}
\usepackage{caption}
\usepackage{pdflscape}
\usepackage{hyperref}
\usepackage{xcolor}
\def\bz{\mathbf{z}}
\def\re{\mathbb{R}}

\def\bbG{\mathbb{G}}
\def\mbi{\mathbb{I}}

\def\bit{\begin{itemize}}
\def\eit{\end{itemize}}

\def\Niin{\mathcal{N}_i^{\mathrm{in}}}
\def\Niout{\mathcal{N}_i^{\mathrm{out}}}

\def\Nikin{\mathcal{N}_{ik}^{\mathrm{in}}}
\def\Nikout{\mathcal{N}_{ik}^{\mathrm{out}}}

\def\1{\mathbf{1}}
\def\0{\mathbf{0}}
\def\T{\mathsf{T}}
\def\la{\langle}
\def\ra{\rangle}

\def\a{{\alpha}}
\def\b{\beta}
\def\g{{\gamma}}

\makeatletter
\newcommand{\eqnum}{\refstepcounter{equation}\textup{\tagform@{\theequation}}}

\newcommand{\argmin}{\operatornamewithlimits{argmin}}

\def\bw{\mathbf{w}}
\def\bg{\mathbf{g}}
\def\bp{\mathbf{p}}

\def\ba{{\mathbf{a}}}

\def\bp{{\mathbf{p}}}
\def\bx{{\mathbf{x}}}
\def\by{{\mathbf{y}}}
\def\bz{{\mathbf{z}}}

\def\bv{{\mathbf{v}}}
\def\bs{{\mathbf{s}}}

\def\D{\mathrm{D}}
\def\K{\mathrm{K}}
\def\G{\mathbb{G}}

\def\diag{{\rm diag}}

\newtheorem{theorem}{Theorem}
\newtheorem{definition}{Definition}
\newtheorem{proposition}{Proposition}
\newtheorem{lemma}{Lemma}
\newtheorem{remark}{Remark}
\newtheorem{assumption}{Assumption}

\def\rev#1{{\color{black}#1}} 

\usepackage{amsmath,amssymb,amsfonts}
\usepackage{algorithmic}
\usepackage{graphicx}
\usepackage{textcomp}

\begin{document}
\title{Accelerated $AB$/Push-Pull Methods for Distributed Optimization over Time-Varying Directed Networks}

\author{Duong Thuy Anh Nguyen, \IEEEmembership{Student Member, IEEE},
Duong Tung Nguyen, \IEEEmembership{Member, IEEE}, \\and Angelia Nedi\'c, \IEEEmembership{Member, IEEE}
\thanks{The authors are with the School of Electrical, Computer and Energy Engineering, Arizona State University, Tempe, AZ, United States. 
Email: \{dtnguy52,~duongnt,~Angelia.Nedich\}@asu.edu. 
This material is based in part upon work supported by the NSF award CCF-2106336.
\textit{Corresponding Author}: Duong Thuy Anh Nguyen.} 
}
\maketitle

\begin{abstract}  
This paper investigates a novel approach for solving the distributed optimization problem in which multiple agents collaborate to find the global decision that minimizes the sum of their individual cost functions. First, the $AB$/Push-Pull gradient-based algorithm is considered, which employs row- and column-stochastic weights simultaneously to track the optimal decision and the gradient of the global cost function, ensuring consensus on the optimal decision. Building on this algorithm, we then develop a general algorithm that incorporates acceleration techniques, such as heavy-ball momentum and Nesterov momentum, as well as their combination with non-identical momentum parameters. Previous literature has established the effectiveness of acceleration methods for various gradient-based distributed algorithms and demonstrated linear convergence for static directed communication networks. In contrast, we focus on time-varying directed communication networks and establish linear convergence of the methods to the optimal solution, when the agents' cost functions are smooth and strongly convex. Additionally, we provide explicit bounds for the step-size value and momentum parameters, based on the properties of the cost functions, the mixing matrices, and the graph connectivity structures. Our numerical results illustrate the benefits of the proposed acceleration techniques on the $AB$/Push-Pull algorithm.
\end{abstract}



\allowdisplaybreaks
\section{Introduction}
    Distributed optimization has attracted significant interest in recent years due to its wide range of applications in large-scale multi-agent systems, such as sensor networks \cite{Rabbat2004}, formation control \cite{Stipanovic2002}, and machine learning \cite{Tsianos2012ML}. In these systems, data samples are distributed across multiple agents with computation tasks divided among them. Communication between agents only occurs through established communication links. This paper considers a system of $n$ agents, where each agent's local cost $f_i$ is determined by its data sample. The goal is for the agents to collaborate and reach a consensus on an optimal solution for the global cost $f$, by solving the following optimization problem:
	\begin{equation} \label{eq-problem}
	\min_{x\in \re^p}~ f(x)=\frac{1}{n}\sum\limits_{i=1}^n f_i(x).
	\end{equation}
    
    The use of decentralized and collaborative approaches for solving the optimization of the sum of convex functions has garnered significant attention in the literature, with many algorithms proposed, including gradient-based methods \cite{Nedic2009,Nedic2011,Ram2009, Srivastava2011,xin2018linear}, dual averaging methods \cite{Duchi2012}, ADMM \cite{Wei2014}, and Newton methods \cite{varagnolo2016newton,mokhtari2017network}. Early works often assume doubly-stochastic weights, which require underlying networks to be undirected or balanced \cite{Shi2015,Wei2014,olshevsky2017linear,pu2018distributed}. To address directed graphs, \cite{Tsianos2012} introduced subgradient-push algorithm, \rev{a decentralized subgradient method based on the push-sum technique,} which is later extended to time-varying graphs in~\cite{nedic2015distributed} with a convergence rate of $O(\ln{t}/\sqrt{t})$ for diminishing step-sizes. Algorithms such as ADD-OPT \cite{xi2018add} and Push-DIGing \cite{nedic2017achieving} further improve the convergence rate by incorporating the push-sum protocol with a gradient estimation approach. These methods require knowledge of agents' out-degree to construct a column-stochastic weight matrix, while algorithms such as \cite{xi2018linear} and FROST \cite{Xin2019FROSTFastRO} only use row-stochastic weights. These algorithms introduce a nonlinear term through division by the Perron eigenvector estimation of the weight matrix, which may result in stability issues. The $AB$/Push-Pull method, introduced by \cite{xin2018linear,pshi21}, eliminates the need for Perron eigenvector estimation by using both row- and column-stochastic weights simultaneously, and demonstrates linear convergence for static directed communication networks. References \cite{Saadatniaki2020} and \cite{Angelia2022AB} further establish linear convergence of this method for time-varying directed graphs, with the latter work also providing an improved analysis and explicit range for the step-size.
    
    The heavy-ball method \cite{POLYAK19641} and Nesterov's momentum \cite{nesterov2003introductory} are popular acceleration techniques for gradient-based methods to achieve faster convergence \cite{Nguyen2023AccGame}. Several distributed algorithms have been proposed in the literature that incorporate these momentum methods. In~\cite{Qu2020}, two variants of an accelerated distributed Nesterov gradient method are proposed for convex (and strongly convex) smooth objective function when the communication network is static undirected. For a static directed network, papers~\cite{Xin2020Heavy} and \cite{Xin2019Nes} propose methods that combine the $AB$/Push-Pull gradient tracking method with a heavy-ball momentum and Nesterov momentum, respectively. In particular, the linear convergence for the 
    $AB$/Push-Pull method with a heavy-ball momentum term is proved in \cite{Xin2020Heavy}, while~\cite{Xin2019Nes} only shows convergence through numerical examples for the $AB$/Push-Pull method with 
    a Nesterov momentum term. Reference~\cite{Huaqing2021} further proposes a double-accelerated method based on $AB$/Push-Pull by incorporating both momentum terms, while ~\cite{Lu2021} proposes to combine Nesterov momentum term with the FROST method. All the aforementioned acceleration methods are studied for a static directed graph. For time-varying communication networks,  \cite{Xiasheng2021} proposes to utilize heavy-ball and Nesterov techniques to accelerate the Push-DIGing algorithm~\cite{nedic2017achieving} that uses column-stochastic matrices only. Furthermore, all prior studies on the double-accelerated method \cite{Huaqing2021,Xiasheng2021} mandate that the parameters for heavy-ball and Nesterov acceleration are identical.
    
    In this paper, we focus on a general network setting where agent communication is given by a sequence of time-varying directed graphs. Building on the $AB$/Push-Pull algorithm, we propose a general formulation that incorporates acceleration techniques such as heavy-ball momentum and Nesterov momentum, as well as their combination with non-identical momentum parameters. This is an innovative departure from existing works, which mandate identical momentum parameters. The $AB$/Push-Pull algorithm does not rely on Perron eigenvector estimation, making it well-suited for time-varying weights and serving as the foundation for the development of these acceleration methods. A key challenge in the analysis is the time-varying nature of the mixing matrices. Our analysis uses time-varying weighted averages and norms to establish consensus contractions for each update step for both row- and column-stochastic mixing matrices. 
    
    We prove that the accelerated algorithm converges linearly to the optimal solution when the agents' cost functions are smooth and strongly convex. The convergence result is derived based on appropriate values for constant step-size and momentum parameters, with explicit upper bounds provided in terms of the cost function properties, mixing matrices, and graph connectivity structures. Numerical results demonstrate that the acceleration of the $AB$/Push-Pull method leads to faster convergence. The results also show that  allowing for different values of acceleration parameters provides increased flexibility and the potential for faster convergence. Our main  contributions can be summarized as follows:
    \begin{itemize}
        \item We propose a novel algorithm that combines the $AB$/Push-Pull method with acceleration techniques such as heavy-ball momentum, Nesterov momentum, and a combination of both, using \textit{non-identical} coefficients.
        \item \rev{We consider a general, \textit{directed}, and \textit{time-varying} communication network} and rigorously prove the linear convergence of the proposed algorithm. The proof extends and improves the previous results \cite{Saadatniaki2020,Angelia2022AB,Xin2020Heavy,Xin2019Nes,Huaqing2021} by providing a comprehensive analysis and explicit ranges for the step-size and momentum parameters in terms of cost function properties and communication structures.
    \end{itemize}
		
    The structure of this paper is as follows. We outline the distributed optimization problem in Section~\ref{sec:Formulation}. In Section~\ref{sec:algo}, we introduce the accelerated algorithm, and its convergence analysis is presented in Section~\ref{sec:result}. The performance of the proposed algorithm is demonstrated in Section~\ref{sec:simulation}. Finally, we conclude with some key points in Section~\ref{sec:conclusion}.

    \textbf{Notations.} Unless otherwise stated, all vectors are considered to be column vectors. The standard Euclidean norm is denoted by $\|\cdot\|$. The notation $\mathbf{1}$ represents a vector with all entries equal to $1$, and $\mathbb{I}$ denotes the identity matrix. The $i$-th entry of a time-varying vector $u_k$ is denoted by $[u_k]_i$. For a vector $v$, we use notation $\min(v)=\min_i v_i$ and $\max(v)=\max_i v_i$. A vector is considered to be stochastic if its entries are nonnegative and sum to $1$.
    
    To denote the $ij$-th entry of a matrix $A_k$, we write $[A_k]_{ij}$. 
    The minimum positive entry of a nonnegative matrix is denoted by $\min^+(A)$. A nonnegative matrix $A\in\mathbb{R}^{n\times n}$ is considered to be row-stochastic if $A\mathbf{1}=\mathbf{1}$, and a nonnegative matrix $B\in\mathbb{R}^{n\times n}$ is considered to be column-stochastic if $\mathbf{1}^{\T} B=\mathbf{1}^{\T}$. Given a positive vector $\ba=(a_1,\ldots,a_n)\in\re^n$, the $\ba$-weighted norm is as follows:
    \[\|\bx\|_\ba=\sqrt{\sum_{i=1}^n a_i\|x_i\|^2} ~~\hbox{for $\bx=(x_1,\ldots,x_n)\in\re^p\!\!\times\!\cdots\!\times\!\re^p$}.\]
    When $\ba=\1$, the norm $\|\bx\|_\ba$ reduces to the Euclidean norm of $\bx$, and we simply write $\|\bx\|$.
    We have the following inequality, 
    \begin{align}\label{eq-norm-anorm}
    &\|\bx\|\le \sqrt{\tfrac{1}{\min(\ba)}}\, \|\bx\|_{\ba}~~\hbox{for $\bx\in\re^p\!\!\times\!\cdots\!\times\!\re^p$, $\ba>\0$}.
    \end{align}
	We let $[n]=\{1,\ldots,n\}$ for an integer $n\ge1$.
	A directed graph $\mathbb{G}=([n],\mathcal{E})$ is specified by the edge set $\mathcal{E}\subseteq [n]\times [n]$ of ordered pairs of nodes. 
Given a directed graph $\mathbb{G}=([n],\mathcal{E})$, the sets $\Niout$ and $\Niin$ denote the out-neighbors and the in-neighbors of a node $i$, i.e.,
	\begin{center}
	$\Niout=\{j\mid (i,j)\in\mathcal{E}\}$ ~and~
	$\Niin=\{j\mid (j,i)\in\mathcal{E}\}$.
	\end{center}
	\rev{When the graph varies over time, we use a subscript $k$ to indicate the time instance. For example, $\mathcal{E}_k$  denotes the edge-set of a graph $\G_k$, $\Nikin$ and $\Nikout$ denote the in-neighbors and the out-neighbors of a node $i$, respectively, at time $k$.}
	
	A directed graph is {\it strongly connected} if there is a directed path from any node to all other nodes in the graph. Given a directed path, the length of the path is the number of edges in the path. For a strongly connected directed graph $\G=([n],\mathcal{E})$, we use the following definitions:
	\begin{definition} [Graph Diameter] \label{def-diam}
	The diameter $\D(\G)$ is the length of the longest path in a collection of all shortest directed paths connecting all ordered pairs of distinct nodes in $\bbG$.
	\end{definition}
	Let $\bp_{jl}$ denote a \textit{shortest directed path from node $j$ to node $l$}, where $j\ne l$. A collection $\mathcal{P}$ of directed paths in $\G$ is a shortest-path graph covering if $\bp_{jl}\in \mathcal{P}$ and $\bp_{lj}\in \mathcal{P}$ for any two
    nodes $j,~l\in [n]$, $j\ne l$. The \textit{utility of the edge} $(j,l)$ with respect to the covering $\mathcal{P}$ is the number of shortest paths in $\mathcal{P}$ that pass through the edge $(j,l)$. 
    Define $\K(\mathcal{P})$ as the maximum edge-utility in $\mathcal{P}$ taken over all edges in the graph, i.e.,
    $\displaystyle \K(\mathcal{P})=\max_{(j,l)\in\mathcal{E}}\sum_{\bp\in\mathcal{P}}\chi_{\{(j,l)\in\bp\}},$
    where $\chi_{\{(j,l)\in\bp\}}$ is the indicator function taking value 1 when $(j,l)\in\bp$ and, otherwise, taking value 0. Denote by $\mathcal{S}(\G)$  the collection of all possible shortest-path coverings of the graph $\G$.
    \begin{definition} [Maximal Edge-Utility] \label{def-edgeut}
    The maximal edge-utility $\K(\G)$ is the maximum value of $\K(\mathcal{P})$ taken over all possible shortest-path coverings $\mathcal{P}\in\mathcal{S}(\G)$,
    i.e., $\K(\G)=\max_{\mathcal{P}\in\mathcal{S}(\G)}\K(\mathcal{P}).$\\
    \end{definition}

\section{Problem Formulation} \label{sec:Formulation}
We consider a system of $n$ agents that are connected by a communication network, with the aim of collaboratively solving the optimization problem in \eqref{eq-problem}, where each function $f_i: \re^p \rightarrow \re$ represents the cost function known only to agent $i$.
The agents aim to find a globally optimal solution by performing local computations and exchanging information with their neighboring agents through a sequence of directed communication networks, represented by a time-varying graph sequence $\{{\bbG_k}\}$.

At each time step $k$, agents communicate over a directed graph $\bbG_k$, and their updates are governed by two non-negative matrices $A_k$ and $B_k$, which {\it align} with the connectivity structure of the graph $\bbG_k$, in the following sense:
\begin{align}\label{eq-alignA}
\!\![A_k]_{ij}\!>0~\forall j\!\in\Nikin\!\cup\{i\}&,\,[A_k]_{ij}\!=0~\forall j\!\not\in\Nikin\!\cup\{i\},\\
\label{eq-alignB}
\!\![B_k]_{ji}\!>0~\forall j\!\in\Nikout\!\cup\{i\}&,\, [B_k]_{ji}\!=0~\forall j\!\not\in\Nikout\!\cup\{i\}.
\end{align}
Moreover, each matrix $A_k$ is row-stochastic and each matrix $B_k$ is column-stochastic
for all $k\ge 0.$ Additionally, we assume that there exist scalars $a>0$ and $b>0$ such that $\min^+(A_k)\ge a$ and $\min^+(B_k)\ge b$ for all $k\ge0$.

We consider the problem under the following assumptions:
\begin{assumption} \label{asm-graphs}For each $k$, the directed graph $\bbG_k=([n],\mathcal{E}_k)$ is  strongly connected.
\end{assumption}

\begin{remark}\label{rem-Cgraphs}
We can relax Assumption~\ref{asm-graphs} by considering a sequence of $C$-strongly connected graphs, i.e., for every $k\ge 0$, there exists an integer $C\ge 1$ such that the graph formed by the edge set $\mathcal{E}^C_k=\bigcup_{i=kC}^{(k+1)C-1}\mathcal{E}_i$ is strongly connected.
\end{remark}

\begin{assumption} \label{asm-functions}
Each $f_i$ is continuously differentiable and has $L$-Lipschitz continuous gradients, i.e., for some $L>0$,
\begin{equation*}
\|\nabla f_i(x)-\nabla f_i(u)\|\le L \|x-u\|,
\quad\hbox{for all $x,u\in \re^p$}.
\end{equation*}
\end{assumption}
	
\begin{assumption}\label{asm-strconv}
The average-sum function $f=\frac{1}{n}\sum_{i=1}^n f_i$ is $\mu$-strongly convex, i.e., for some $\mu>0$, 
\[\la \nabla f(x)-\nabla f(u),x-u\ra \ge \mu\|x-u\|^2\quad\hbox{for all $x,u\in \re^p$}.\]
\end{assumption}
\begin{remark}
The strong convexity condition implies that problem~\eqref{eq-problem} has a unique optimal solution denoted by $x^*$, i.e.,
\[x^*=\argmin_{x\in\re^p} f(x).\]
\end{remark}

\section{Accelerated $AB$/Push-Pull Methods} \label{sec:algo}
\rev{The $AB$/Push-Pull method, initially proposed in \cite{pshi21,xin2018linear}, is a general framework that unifies many of the existing decentralized first-order methods with gradient tracking. Thus, it is worthwhile to consider adding momentum terms to accelerate its convergence.} In this section, we introduce a comprehensive approach that encompasses three methods for accelerating the distributed $AB$/Push-Pull algorithm for time-varying directed graphs. These methods include the  heavy-ball method~\cite{POLYAK19641} (also, in Section 3.2.1 of~\cite{Polyak}), the Nesterov gradient method~\cite{nesterov2003introductory}, and a combination of these two methods.

\rev{Let each agent $i\in\{1,2,\ldots,n\}$ possess two local copies $x_k^i\in \re^p$ and $s_k^i\in \re^p$ of the decision variable and a gradient-tracking variable $y_k^i\in\mathbb{R}^p$ which is an estimate of a ``global update direction", at iteration $k$. These variables are maintained and updated over time, as follows:} 

\rev{To begin, agents are directed to use the step-size $\alpha>0$, the heavy-ball momentum parameter $\beta \geq 0$, and the Nesterov momentum parameter $\gamma \geq 0$. Additionally, each agent $i$ initializes their updates with arbitrary vectors $x_{-1}^i,x_0^i,s_0^i$ and with $y_0^i=\nabla f_i(s_0^i)$, without the need for coordination among agents. Each agent $i$ also independently decide on the entries of $A_k$ in the $i$-th row for their in-neighbors $j\in\Nikin$, while the value $[B_k]_{ij}$ is determined by agent $j\in \Nikin$. At every time $k$, every agent $i$ sends its vector $s_k^i$ and a scaled direction $[B_k]_{\ell i}y_k^i$ to its out-neighbors $\ell\in\Nikout$. Every agent $i$ also receives these vectors sent by its in-neighbors $j\in\Nikin$.
Upon the information exchange step, every agent $i$ updates its vectors using equations \eqref{eq-x}--\eqref{eq-y} for all $k\ge0$. The proposed procedure is outlined in Algorithm 1.}

\begin{table}[t!]
\vspace{0.2cm}
\centering \normalsize
    \begin{tabular}{l}
    \hline
    \multicolumn{1}{c}{\rev{\textbf{Algorithm 1:} Accelerated $AB$/Push-Pull}}\\
    \hline
    Agents are directed to use $\a>0$, $\b\ge 0$ and $\g\ge 0$.\\
    Every agent $i\in[n]$ initializes with arbitrary initial vectors\\ $x_{-1}^i,x_0^i,s_0^i\in\mathbb{R}^p$ and $y_0^i=\nabla f_i(s_0^i)$.\\
    \textbf{for} $k=0,1,\ldots,$ every agent $i\in[n]$ does the following:\!\!\\
    \emph{ } In-bounds mixing weights $[A_k]_{ij}$, for all $j\in\Nikin$;\\
    \emph{ } Out-bounds pushing weights $[B_k]_{\ell i}$, for all $\ell\in\Nikout$;\\
    \emph{  } Receives $s_k^j$ and $[B_k]_{ij}y_k^j$ from in-neighbors $j\in\Nikin$;\\
    \emph{  } Sends $s_k^i$ and $[B_k]_{\ell i}y_k^i$ to out-neighbors $\ell\in\Nikout$;\\
    \emph{ } Updates $x^i_{k+1}$, $s_{k+1}^i$ and $y_{k+1}^i$ by \\
    \emph{~~~~} $\displaystyle x_{k+1}^i =  \sum_{j=1}^n[A_k]_{ij}s_k^j  - \a y_k^i + \b(x_k^i-x_{k-1}^i),$ ~\quad\eqnum \label{eq-x}\\
    \emph{~~~~} $\displaystyle s_{k+1}^i =  x_{k+1}^i + \g(x_{k+1}^i-x_k^i),$ \qquad\qquad\qquad\quad\,\eqnum\label{eq-s}\\
    \emph{~~~~} $\displaystyle y_{k+1}^i =  \sum_{j=1}^n [B_k]_{ij}y_k^j + \nabla f_i(s_{k+1}^i) -\nabla f_i(s_k^i),$ \quad\eqnum\label{eq-y}\\
    \textbf{end for}\\
    \hline
    \end{tabular}
    \vspace{-0.5cm}
\end{table}

The method in Algorithm 1 is a generalization of three methods to accelerate the $AB$/Push-Pull algorithm, namely,
\begin{itemize}
    \item[\tiny$\bullet$] $\b>0,\g=0$: Heavy-ball method
    \item[\tiny$\bullet$] $\b=0,\g>0$: Nesterov gradient method
    \item[\tiny$\bullet$] $\b>0,\g>0$: Combination of Nesterov gradient and heavy-ball methods.
\end{itemize}

\noindent \rev{\textbf{Relations to the $AB$/Push-Pull algorithm:} From the viewpoint of an agent, the information about the gradients is pushed to the neighbors, while the information about the decision variable is pulled from the neighbors, as noted in \cite{pshi21}.} Hence, the update step for aggregating the decision variables using the row-stochastic matrix $A_k$ is referred to as a {\it pull-step}, while the step for tracking the average gradients using the column-stochastic matrix $B_k$ is referred to as a {\it push-step} as it mimics the push-sum consensus method, originally proposed in~\cite{Kempe2003}, later used in decentralized gradient-based methods~\cite{Tsianos2012,Tsianos2012ML,nedic2015distributed}, and recently studied in~\cite{Rezaienia2020}. \rev{Moreover, the role of gradient tracking is to account for the heterogeneity of the local data distributions among agents.}

\noindent \rev{\textbf{Acceleration techniques:}} Intuitively, the heavy-ball momentum term, represented by $\b(x_k^i-x_{k-1}^i)$, accelerates the gradient method by adding inertia to the updates. The Nesterov momentum term, represented by $\g(x_{k+1}^i-x_k^i)$, is mathematically shown to improve convergence rate, but its underlying intuition is not immediately clear. A geometric interpretation of the Nesterov accelerated algorithm can be found in \cite{bubeck2015geometric}. We discuss the impact of incorporating the momentum terms on the convergence rate of the standard gradient method, i.e.,
\[x_{k+1}=x_{k}-\a \nabla f(x_{k}),\]
where $\a>0$ is a stepsize. We recall the following updates of the heavy-ball method~(see\cite{POLYAK19641}, or Section 3.2.1 of \cite{Polyak}):
\[x_{k+1}=x_{k}-\a \nabla f(x_{k})+\b(x_{k}-x_{k-1}),\]
where $\b>0$,
and the Nesterov gradient method~\cite{nesterov2003introductory}:
\begin{align*}
    x_{k+1}&=y_{k}-\a \nabla f(y_{k}),\\
    y_{k+1}&=x_{k+1}+\g(x_{k+1}-x_{k}),
\end{align*}
where $\gamma>0$.
The convergence rate of the gradient method is well-known to be $\mathcal{O}\big(\big(\tfrac{Q-1}{Q+1}\big)^k\big)$, where $Q=\tfrac{L}{\mu}$ is the condition number of the objective function $f(x)$ (see~\cite{Polyak}). By properly choosing the parameters $\a$ and $\b$, the heavy-ball method can achieve a faster, locally accelerated rate of $\mathcal{O}\big(\big(\tfrac{\sqrt{Q}-1}{\sqrt{Q}+1}\big)^k\big)$. The Nesterov gradient method also has the rate of $\mathcal{O}\big(\big(\tfrac{\sqrt{Q}-1}{\sqrt{Q}}\big)^k\big)$ when $0< \a \le 1/L$ and $\g = (\sqrt{Q}-1)/(\sqrt{Q}+1)$, which is faster than the gradient method but slower than the heavy-ball method in terms of the dependence on the condition number~$Q$ (since $Q\ge 1$).
These acceleration methods, as well as their combination, have been studied in the context of distributed optimization algorithms. It has been shown numerically and/or theoretically in \cite{Xin2020Heavy,Xin2019Nes,Huaqing2021} that the accelerated  $AB$/Push-Pull algorithms converge linearly for static directed graphs with appropriate choices of step-size and momentum parameters. However, despite this progress, a comprehensive theoretical analysis for the convergence on time-varying directed graphs is still an open problem, which will be addressed in the next section.\\



\section{Convergence Analysis}\label{sec:result}
In this section, we present a detailed convergence analysis of the generic accelerated $AB$/Push-Pull algorithm over a time-varying directed communication network. We begin by introducing some preliminary results. We then proceed to establish the estimates for the four quantities of interest, namely the optimality gap, the consensus error, the state difference, and the gradient estimation error. By combining these estimates, we demonstrate the linear convergence of the algorithm. Finally, we provide a step-size selection rule and the bounds for the acceleration parameters that ensure the convergence of the proposed algorithm.

\subsection{Preliminary Results}
We first present the contraction property of the gradient mapping, assuming that the objective function is strongly convex and has Lipschitz continuous gradients, in the following:
\begin{lemma}[\cite{Polyak}]\label{lem-contraction} 
Let $f$ be a $\mu$-strongly convex and $L$-smooth function. For $0<\a<2L^{-1}$, we have
\[\|x-x^*-\a \nabla f(x)\|\le q(\a) \|x-x^*\|\qquad\hbox{for all $x$},\]
where $q(\a)=\max\{|1-\a\mu|,|1-\a L|\}<1$.
\end{lemma}

We then proceed by presenting some foundational results that will support our later analysis.
\begin{lemma}[\cite{nguyen2022distributed}, Corollary~5.2]\label{lem-normlincomb}
Consider a vector collection $\{u_i, \, i\in[n]\}\subset\re^p$, and a scalar collection $\{\g_i,\, i\in[n]\}\subset\re$ of scalars such that $\sum_{i=1}^n \g_i=1$. For all $u\in \re^p$, we have the following relation:
\begin{align*}
\Bigg\|\!\sum_{i=1}^n \g_i u_i - u \Bigg\|^2 \!\!= \sum_{i=1}^n \g_i \|u_i-u\|^2 -\!\sum_{i=1}^n \g_i \Bigg\|u_i \!-\! \sum_{\ell=1}^n \g_\ell u_\ell\Bigg\|^2\!\!.
\end{align*}
\end{lemma}

\begin{lemma}[\cite{nguyen2022distributed}, Lemma 5.4] \label{lem-amatrices}
Consider a sequence $\{A_k\}$ of row-stochastic matrices. Then, there exists a sequence $\{\phi_k\}$ of stochastic vectors such that
\begin{align}\label{eq-phik}
\phi_{k+1}^{\T}A_k=\phi_k^{\T}\qquad\hbox{for all $k\ge0$}.
\end{align}
Moreover, if Assumption~\ref{asm-graphs} holds
and $A_k$ is aligned with the graph $\mathbb{G}_k$ (see~\eqref{eq-alignA}) with  $\min^+(A_k)\ge a>0$ for all $k\ge0$, then
the entries of each $\phi_k$ have a uniform lower bound, i.e.,	$[\phi_k]_i\ge \tfrac{a^n}{n}$ for all $i\in[n]$ and for all $k\ge0$.
\end{lemma}

\begin{lemma}[\cite{Angelia2022AB}, Lemma 3.4] \label{lem-bmatrices}
Consider a sequence $\{B_k\}$ of column-stochastic matrices and the vector sequence  $\{\pi_k\}$, defined as follows:
\begin{align}\label{eq-pik}
\pi_{k+1}=B_k\pi_k,\qquad\hbox{initialized with }\ \pi_0=\tfrac{1}{n}\1.
\end{align} 
Then, the vectors $\pi_k$ are stochastic vectors. Moreover, if Assumption~\ref{asm-graphs} holds, 
where $B_k$ is aligned with the graph $\mathbb{G}_k$  and  $\min^+(B_k)\ge b>0$ for all $k\ge0$, then
$[\pi_k]_i\ge \tfrac{b^n}{n}$ for all $i\in[n]$ and $k\ge0$.
\end{lemma}

\begin{remark} \label{rem-CgraphsAB}
    \rev{When the graph sequence is $C$-strongly-connected (Remark~\ref{rem-Cgraphs}), the product of weight matrices $A_{k+C-1}\ldots A_{k+1}A_{k}$ and  $B_{k+C-1}\ldots B_{k+1}B_{k}$ are row- and column-stochastic, respectively. These matrices represent the directed paths among the nodes in the composition of the graphs $\G_{k},\ldots,\G_{k+C-1}$, capturing the underlying connectivity patterns and facilitating the understanding of information flow dynamics.} Moreover, the more general results of Lemma~\ref{lem-amatrices} and Lemma~\ref{lem-bmatrices} indicate the existence of stochastic vector sequences ${\phi_k}$ and ${\pi_k}$, such that for all $k \ge 0$, 
    \begin{align*}
    &\phi'_{k+C}\left(A_{k+C-1}\ldots A_{k+1}A_{k}\right)=\phi'_k \\
    \text{and}~ &\pi_{k+C}=\left(B_{k+C-1}\ldots B_{k+1}B_{k}\right)\pi_k.
    \end{align*}
  Furthermore,
    $$[\phi_k]_i\ge\dfrac{a^{nC}}{n} ~\text{and}~[\pi_k]_i\ge\dfrac{b^{nC}}{n} \quad \text{for all } i\in[n].$$
\end{remark}

Consider a strongly connected directed graph $\bbG=([n],{\cal E})$, and weight matrices $A$ and $B$ that are aligned with the graph $\bbG$ (in sense of equations~\eqref{eq-alignA} and~\eqref{eq-alignB}). Let $\mathsf{D}(\bbG)$ and $\mathsf{K}(\bbG)$ be the diameter and the maximal edge-utility of the graph $\bbG$, respectively. We have the following two results:
\begin{lemma}[\cite{nguyen2022distributed}, Lemma 6.1]\label{lem-basic-xcontract}
Let $A$ be a row-stochastic matrix, $\phi$  be a stochastic vector and let $\pi$ be a nonnegative vector such that $\pi^{\T}A=\phi^{\T}$. Consider a collection of vectors $x_1,\ldots,x_n\in\re^p$. For $\hat x_\phi=\sum_{i=1}^n \phi_i x_i$, we have
\begin{align*}
\sqrt{\sum_{i=1}^n \pi_i\Bigg\|\sum_{j=1}^n A_{ij} x_j- \hat x_\phi\Bigg\|^2} &\le c\sqrt{\sum_{j=1}^n \phi_j \|x_j - \hat x_\phi \|^2},
\end{align*}
where $c=\sqrt{1 -   \frac{\min(\pi) (\min^+(A))^2}{\max^2(\phi)\mathsf{D}(\bbG)\mathsf{K}(\bbG)}}\in (0,1)$ is a scalar.
\end{lemma}

\begin{lemma}[\cite{Angelia2022AB}, Lemma 4.5]\label{lem-basic-ycontract}
Let $B$ be a column-stochastic matrix, $\nu$  be a stochastic vector with positive entries, i.e., $\nu_i>0,~\forall i\in[n]$, and let the vector $\pi$ be given by $\pi=B\nu$. Consider the vectors $y_1,\ldots,y_n\in\re^p$ and vectors $w_i=\sum_{j=1}^n B_{ij} y_j$ for all $i\in[n]$, we have
\begin{align*}
	\sqrt{\sum_{i=1}^n\pi_i \Bigg\|\frac{w_i}{\pi_i} - \sum_{\ell=1}^m y_\ell\Bigg\|^2}
	\le\tau\,
	\sqrt{\sum_{i=1}^n \nu_i \Bigg\|\frac{y_i}{\nu_i} - \sum_{\ell=1}^n y_\ell\Bigg\|^2},
\end{align*}	
where $\tau=\sqrt{1 -   \frac{\min^2(\nu)\,(\min^+(B))^2}{\max^2(\nu) \max(\pi)\, \mathsf{D}(\bbG)\mathsf{K}(\bbG)}}\in (0,1)$ is a scalar.
\end{lemma}

The column stochastic property of the matrices $B_k$ ensures that the sum of the $y$-iterates, at any time $k$, is equal to the sum of the gradients $\nabla f_i(s_k^i)$, as seen in the following lemma (the proof follows from mathematical induction on $k$, and uses the column-stochasticity of $B_k$ and the initialization of the $y$-variables; see~Algorithm 1). 
\begin{lemma} \label{lem-sumgrad}
Consider Algorithm 1, and assume that each $B_k$ is column-stochastic. Then, we have
\[\sum_{i=1}^n y_k^i=\sum_{i=1}^n \nabla f_i(s_k^i)\quad\hbox{for all $k\ge0$.}\]
\end{lemma}

\subsection{Main Results}
The convergence analysis of the accelerated $AB$/Push-Pull method is based on establishing a contraction relationship between the following four quantities:
(i) the optimality gap, (ii) the consensus error, (iii) the state difference, and (iv) the gradient estimation error, given respectively as follows:
\begin{subequations}\label{eq-quants}
\begin{align}
&\!\!\!\|\hat x_k -x^*\|, ~~ D(\bx_k,\phi_k)= \sqrt{\sum_{i=1}^n[\phi_k]_i\|x_k^i-\hat x_k\|^2},\label{eq-x-D-quants}\\
&\!\!\!\|\bx_{k}\!-\!\bx_{k-1}\|,~ S(\by_k,\!\pi_k) \!=\!\! \sqrt{\sum_{j=1}^n[\pi_k]_j\Bigg\|\frac{y_k^j}{[\pi_k]_j} \!- \!\!\sum_{\ell=1}^n y_k^\ell\Bigg\|^2}\!,\!\! \label{eq-x-S-quants}
\end{align}
\end{subequations}
where $\hat x_k= \sum_{i=1}^n[\phi_k]_i x_k^i$, $\bx_k=(x_k^1,\ldots,x_k^n)$ and $x^*$ is the optimal solution of problem~\eqref{eq-problem}.  We define the constants $\varphi_k>0$, $r_k>0$, $c_k\in(0,1)$ and $\tau_k\in(0,1)$ as follows
\begin{align}\label{eq-const-all}
&r_k=\sqrt{n} + \tfrac{1}{\sqrt{\min(\pi_{k+1})}},~
c_k=\sqrt{1- \tfrac{\min(\phi_{k+1})\, a^2} {\max^2(\phi_k)\,\mathsf{D}(\bbG_k)\mathsf{K}(\bbG_k)}}, \nonumber\\
&\varphi_k=\sqrt{\tfrac{1}{\min(\phi_{k})}},~ \tau_k =\sqrt{1 -   \tfrac{\min^2(\pi_k)\,b^2}
{\max^2(\pi_k) \max(\pi_{k+1}) \mathsf{D}(\bbG_k)\mathsf{K}(\bbG_k)}}.
\end{align}

We first establish the recursive relation for the weighted average $\{\hat x_k\}$ that will be utilized in the subsequent analysis.
\begin{lemma} \label{lem-weightedavg}
The weighted average sequence $\{\hat x_k\}$ satisfies,
\begin{align*}
	\hat x_{k+1}&=\hat x_{k} +\sum_{j=1}^n (\b[\phi_{k+1}]_j+\g[\phi_{k}]_j)(x_k^i-x_{k-1}^i) \nonumber\\
	&- \a \sum_{i=1}^n [\phi_{k+1}]_i y_k^i,\qquad\hbox{for all $k\ge0$.}
\end{align*}
\end{lemma}
\begin{proof}
	Plugging in the update for $s_k^j$ from \eqref{eq-s} into the update for $x_k^j$ in \eqref{eq-x} and rearranging the terms results in:
	\begin{align}\label{eq-xcomb}
	\!\!\!x_{k+1}^i \!=\!\! \sum_{j=1}^n[A_k]_{ij}x_k^i\!-\! \a y_k^i \!+\! \Bigg(\!\!\b+\g\!\sum_{j=1}^n[A_k]_{ij}\!\!\Bigg)(x_k^i\!\!-\!x_{k-\!1}^i).\!\!\!
	\end{align}
	By taking a weighted average of $x_{k+1}^i$ using the $\phi_{k+1}$ weights, from relation~\eqref{eq-xcomb} we obtain
	\begin{align*}
	\sum_{i=1}^n [\phi_{k+1}]_i x_{k+1}^i&=
 \sum_{i=1}^n [\phi_{k+1}]_i\sum_{j=1}^n[A_k]_{ij}x_k^j - \a \sum_{i=1}^n [\phi_{k+1}]_i y_k^i\\
	&+\sum_{i=1}^n[\phi_{k+1}]_i (\b+\g\sum_{j=1}^n[A_k]_{ij})(x_k^i-x_{k-1}^i).
	\end{align*}
	Since $\phi_{k+1}^{\T}A_k=\phi_k^{\T}$, for the double-sum term we have
	\[\!\!\sum_{i=1}^n [\phi_{k\!+\!1}]_i \!\!\sum_{j=1}^n [A_k]_{ij}x_k^j
	\!=\!\!\!\sum_{j=1}^n \!\!\Big(\sum_{i=1}^n [\phi_{k\!+\!1}]_i [A_k]_{ij}\!\!\Big) x_k^j\!\!
	=\!\!\!\sum_{j=1}^n [\phi_k]_j x_k^j.\]
	By using the definition of $\displaystyle\hat x_{k}\!=\!\!\sum_{j=1}^n [\phi_k]_j x_k^j$, and combining the preceding two equations, we arrive at the desired relation.
\end{proof}

We next examine the behavior of the directions $y_k^i$ generated by the update in~\eqref{eq-y}. The analysis will make use of some weighted norms of scaled directions $y_k^i$, where the scalings are time-varying and defined by a stochastic vector sequence $\{\pi_k\}$ associated with the matrix sequence $\{B_k\}$. 
\begin{lemma} \label{lem-yinv}
Under Assumption~\ref{asm-graphs},
we have,
\begin{equation*}
	\sqrt{\sum_{i=1}^n \frac{\|y_k^i\|^2}{[\pi_k]_i}}
	\le S(\by_k,\pi_k) +\Bigg\|\sum_{\ell=1}^n y_k^\ell\Bigg\|\qquad\hbox{for all $ k\ge 0$}.
\end{equation*}
\end{lemma}
\begin{proof}
    By Lemma~\ref{lem-bmatrices}, we have $\pi_k > 0$ for all $k\ge0$, thus
	\[\sum_{i=1}^n \frac{\|y_k^i\|^2}{[\pi_k]_i} 
	=\sum_{i=1}^n [\pi_k]_i\,\left\|\frac{y_k^i}{[\pi_k]_i} \right\|^2.\]
	Using the relation in Lemma~\ref{lem-normlincomb} with $\g_i=[\pi_k]_i$ and
	$u_i=y_k^i/[\pi_k]_i$ for all $i$, and $u=0$, we obtain
	\begin{align*}
	\sum_{i=1}^n [\pi_k]_i\,\Bigg\|\frac{y_k^i}{[\pi_k]_i} \Bigg\|^2
	&=\sum_{i=1}^n [\pi_k]_i\Bigg\|\frac{y_k^i}{[\pi_k]_i} -\sum_{\ell=1}^n y_k^\ell\Bigg\|^2 
	+\Bigg\|\sum_{\ell=1}^n y_k^\ell\Bigg\|^2\\
	&=S^2(\by_k,\pi_k) +\Bigg\|\sum_{\ell=1}^n y_k^\ell\Bigg\|^2,
	\end{align*}
	where the last equality is implied from the definition of the $S$-quantity in~\eqref{eq-x-S-quants}. Therefore, 
	\begin{equation*}
	\sqrt{\sum_{i=1}^n \!\frac{\|y_k^i\|^2}{[\pi_k]_i}}\!	\!=\!\!\sqrt{\!S^2(\by_k,\!\pi_k) \!+\!\Bigg\|\!\sum_{\ell=1}^n y_k^\ell\Bigg\|^2}
	\!\!\le \!S(\by_k,\!\pi_k) +\Bigg\|\sum_{\ell=1}^n y_k^\ell\Bigg\|,\end{equation*}
	where the inequality is obtained by using $\sqrt{a+b}\le \sqrt{a}+\sqrt{b}$, 
	which is valid for any two scalars
	$a,b\ge0$.
\end{proof}

\begin{lemma} \label{lem-ysum}
Under Assumption~\ref{asm-graphs}, Assumption~\ref{asm-functions}, and Assumption~\ref{asm-strconv},
we have the following relation for the sum of the $y$-iterates, which holds for all $k \geq 0$:
\[\!\!\left\|\sum_{i=1}^n y_k^i\right\|
\!\!\le \!\!L\sqrt{n} \left(\varphi_k\|\hat x_k \!-\! x^*\| \!\!+\!\!\varphi_k D(\bx_k,\!\phi_k)\!\!+\!\g\|\bx_k \!\!-\!\bx_{k-1}\|\right)\!,\]
where $\varphi_k>0$ is as given in~\eqref{eq-const-all}.
\end{lemma}

\begin{proof}
    By Lemma~\ref{lem-sumgrad}, we have 
	\[\Bigg\| \sum_{i=1}^n y_k^i\Bigg\|=\Bigg\|\sum_{i=1}^n \nabla f_i(s_k^i)\Bigg\|
	=\Bigg\|\sum_{i=1}^n\Big(\nabla f_i(s_k^i) -\nabla f_i(x^*)\Big) \Bigg\|,\]
where we use $\sum_{i=1}^n\nabla f_i(x^*)=0$, valid for the optimal solution $x^*$ to problem~\eqref{eq-problem} (which exists and  is unique due to Assumption~\ref{asm-strconv}). By Assumption~\ref{asm-functions} that each $f_i$ has Lipschitz continuous gradients with a constant $L>0$, we obtain
	\[\Bigg\| \sum_{i=1}^n y_k^i\Bigg\|
	\le \sum_{i=1}^n\Bigg\| \nabla f_i(s_k^i)-\nabla f_i(x^*)\Bigg\|
	\le L\sum_{i=1}^n\| s_k^i -x^*\|.\]
	We define $\bs_k=(s_k^1,\ldots,s_k^n)$ and $\bx^*=(x^*,\ldots,x^*)$. By H\"older's inequality we have $\sum_{i=1}^n a_i\le\sqrt{n\sum_{i=1}^n a_i^2}$, for all $a_i$, $i\in[n]$, implying that
	\begin{align}\label{eq-sumy1}
	\Bigg\| \sum_{i=1}^n y_k^i\Bigg\| &\!\le\! L\sqrt{n}\|\bs_k \!-\!\bx^*\|
	\!\le\! L\sqrt{n}(\|\bx_k \!-\!\bx^*\|\!+\!\g\|\bx_k \!-\!\bx_{k-1}\|)\nonumber\\
	&\!\le\! L\sqrt{n} (\varphi_k|\bx_k -\bx^*\|_{\phi_k}+\g\|\bx_k -\bx_{k-1}\|),\!\!
	\end{align}
where the second inequality follows from the updates of the $s$-iterates in~\eqref{eq-s}, while the last inequality follows from the relation for the norms in~\eqref{eq-norm-anorm} and the fact that $[\phi_k]_i>0$ for all $i$ and $k$ (by Assumption~\ref{asm-graphs} and  Lemma~\ref{lem-amatrices}).
	For the quantity $\|\bx_k -\bx^*\|_{\phi_k}$, applying the relation in Lemma~\ref{lem-normlincomb} with $u_i=x_k^i$, $\g_i=[\phi_k]_i$ for all $i$, and $u=x^*$ yields
	\[\sum_{i=1}^n [\phi_k]_i\|x_k^i-x^*\|^2
	=\|\hat x_k-x^*\|^2 + \sum_{i=1}^n [\phi_k]_i\|x_k^i-\hat x_k\|^2,\]
	where $\hat x_k=\sum_{\ell=1}^n[\phi_k]_\ell x_k^\ell$. Using the definition of $D(\bx_k,\phi_k)$ in~\eqref{eq-quants}, we further derive the following inequality:
	\begin{equation}\label{eq-xrel2}
	\left\|\bx_k - \bx^*\right\|_{\phi_k}
	\le  \|\hat x_k-x^*\|+D(\bx_k,\phi_k),
	\end{equation}
	which uses $\sqrt{a+b}\le \sqrt{a}+\sqrt{b}$. The desired relation follows by combining the relations in \eqref{eq-sumy1} and~\eqref{eq-xrel2} with the inequality $\|\bx_k -\bx_{k-1}\|_{\phi_k} \le \|\bx_k -\bx_{k-1}\|$.
\end{proof}

In the following, we derive upper bounds for each of the four quantities defined in \eqref{eq-quants}. We begin by assessing the optimality gap in the subsequent proposition. 

\begin{proposition}\label{prop-waverx}
Let Assumption~\ref{asm-graphs}, Assumption~\ref{asm-functions}, and Assumption~\ref{asm-strconv} hold.
Let the step-size $\a$ in Algorithm 1 be such that $0<\a<\tfrac{2}{n\min(\pi_k)L}$, where $L$ is the gradient Lipschitz constant. Then, we have for all $k\ge0$: 
\begin{align*}
	\|\hat x_{k+1}-&x^*\| \le  q_k(\a) \|\hat x_k - x^*\| + \a L\sqrt{n}\varphi_kD(\bx_k,\phi_k)\\
	&+\a S(\by_k,\pi_k)+\left(\b+(1+\a L\sqrt{n})\g\right)\|\bx_{k}-\bx_{k-1}\|,
	\end{align*}
	where $q_k(\a) \!=\! \max\{|1\!-\a n \min(\pi_k)\mu|,|1\!-\a n \min(\pi_k)L|\}$.
\end{proposition}
\begin{proof}
See Appendix~\ref{app:Prop1}.
\end{proof}

In the next proposition, we investigate the behavior of the deviation of the iterates $x_k^i$, $i\in[n]$, from their weighted average $\hat x_k$, as measured by the $\phi_k$-weighted dispersion $D(\bx_k,\phi_k)$.
\begin{proposition}\label{prop-xcontract} 
Under Assumption~\ref{asm-graphs}, 
Assumption~\ref{asm-functions}, and Assumption~\ref{asm-strconv}, the following inequality holds for all $ k\ge 0$,
\begin{align*}
	D(\bx_{k+1},\phi_{k+1})\!\le \! 
	 (c_k+\a L\sqrt{n}\varphi_k) D(\bx_k,\phi_k)\!+\!\a S(\by_k,\pi_k)\\
	+\a L\sqrt{n}\varphi_k\|\hat x_k \!-\! x^*\| \!+\! (\b\!+\!\g (c_k+ \a L\sqrt{n}))\|\bx_{k}\!-\!\bx_{k-1}\|.
\end{align*}
\end{proposition}
\begin{proof}
See Appendix~\ref{app:Prop2}.
\end{proof}

The next result establishes an upper bound for the state difference of the $x$-sequence produced by the update in~\eqref{eq-x}.
\begin{proposition}\label{prop-xdiff}
Let Assumption~\ref{asm-graphs},
Assumption~\ref{asm-functions}, and Assumption~\ref{asm-strconv} hold. Then, for all $ k\ge 0$, we have:
\begin{align*}
&\|\bx_{k+1}\!\!-\!\bx_{k}\| \!\le\! (\b\!+\!\g \sqrt{n}(1\!+\!\a L))\|\bx_{k}\!\!-\!\bx_{k-1}\|\!+\!\a S(\by_k,\!\pi_k)\\
&\!+\!\a L\sqrt{n}\varphi_k\|\hat x_k \!-\! x^*\|\!+\!(c_k\varphi_{k+1}\!\!+\!\varphi_k\!+\!\a L\sqrt{n}\varphi_k) D(\bx_k,\!\phi_k).
\end{align*}
\end{proposition}
\begin{proof} 
    Let us denote $z_k^i=\sum_{j=1}^n [A_k]_{ij} x_k^j$ and $v_k^i=\sum_{j=1}^n [A_k]_{ij} (x_k^j- x_{k-1}^j)$. Define the vectors $\bz_k=(z_k^1,\ldots,z_k^n)$, $\bv_k=(v_k^1,\ldots,v_k^n)$, $\by_k=(y_k^1,\ldots,y_k^n)$ and $\hat \bx_k=(\hat x_k,\ldots,\hat x_k)$. Then, we can write the $x$-update in \eqref{eq-xcomb} compactly as follows:
    \begin{equation}\label{eq-x-compact}
    \bx_{k+1}=\bz_k+\g\bv_k-\a\by_k+\b(\bx_{k}-\bx_{k-1}).\end{equation}
    By adding and subtracting $\hat \bx_{k}$ and using the triangle inequality, we obtain:
	\begin{align}\label{eq-xdiffp1}
	&\|\bx_{k+1}-\bx_k\| \le \|\bx_{k+1} - {\hat \bx}_k\| +\|{\hat \bx}_k-\bx_k\|\nonumber\\
	\!\!\!\le &\|\bz_k\!\!-\!{\hat \bx}_k\|\!\!+\!\a\|\by_k\|\!\!+\!\b\|\bx_{k}\!-\!\!\bx_{k-\!1}\|\!\!+\!\g\|\bv_{k}\|\!+\!\|\bx_k \!\!-\! {\hat \bx}_k\|,\!\!\!\!
	\end{align}
	where the last inequality follows from the compact representation of the $x$-update in~\eqref{eq-x-compact}. 
	For the first term in \eqref{eq-xdiffp1}, we use the relation for norms in~\eqref{eq-norm-anorm} and Lemma~\ref{lem-basic-xcontract} with $A=A_k$, $x_i=x_k^i$, and $\phi_{k+1}^{\T}A_k=\phi_k^{\T}$ to obtain the following bound:
	\[\|\bz_k-{\hat \bx}_k\|\le \varphi_{k+1}\left\|\bz_k-{\hat \bx}_k\right\|_{\phi_{k+1}} \le c_k\varphi_{k+1}\|\bx_k - {\hat \bx}_k\|_{\phi_k}.\]
	For the second term in \eqref{eq-xdiffp1}, using the fact that the vector $\pi_k$ is stochastic, Lemma~\ref{lem-yinv} and Lemma~\ref{lem-ysum}, we derive
	\begin{align*}
	\|\by_k\| &= \sqrt{\sum_{i=1}^n [\pi_k]_i \frac{\|y_k^i\|^2}{[\pi_k]_i}} \le \sqrt{\sum_{i=1}^n  \frac{\|y_k^i\|^2}{[\pi_k]_i}}\le  S(\by_k,\!\pi_k)\\
	 &+ L\sqrt{n} \left(\varphi_k\|\hat x_k \!-\! x^*\| \!+\! \varphi_kD(\bx_k,\!\phi_k)\!+\!\g\|\bx_k \!-\!\bx_{k-1}\|\right).
    \end{align*}
	For the fourth term in \eqref{eq-xdiffp1}, we have
	\begin{align*}
	\|\bv_{k}\|^2 &= \sum_{i=1}^n \Bigg\|\sum_{j=1}^n [A_k]_{ij} (x_k^j- x_{k-1}^j)\Bigg\|^2\\
	& \le \sum_{i=1}^n \sum_{j=1}^n [A_k]_{ij}\|x_k^j- x_{k-1}^j\|^2 \le n\|\bx_{k}-\bx_{k-1}\|^2,
    \end{align*}
    where we use the fact that $A_k$ is row-stochastic.
    
	The last term in \eqref{eq-xdiffp1} follows from the relation for norms in~\eqref{eq-norm-anorm} and the definition of $D(\bx_k,\phi_k)$, as follows
	\begin{align*}
	\|\bx_k - {\hat \bx}_k\| \le \varphi_{k} \|\bx_k - {\hat \bx}_k\|_{\phi_k} = \varphi_{k}D(\bx_k,\phi_k).
	\end{align*}
	Combining the estimates for each of the  terms in~\eqref{eq-xdiffp1}, we obtain the desired relation.
\end{proof}

Next, we provide a recursive relation for the gradient estimation error $S(\by_k,\pi_k)$, as given in the following proposition.
\begin{proposition}\label{prop-ycontract}
Under Assumption~\ref{asm-graphs} and Assumption~\ref{asm-functions}, the following inequality holds for all $ k\ge 0$,
\begin{align*}
	S(\by_{k+1},\pi_{k+1})
	\le &~\tau_k\, S(\by_k,\pi_k)+L r_k(1+\g)\|\bx_{k+1}-\bx_k\|\\
	&+L r_k\g\|\bx_{k}-\bx_{k-1}\|.
\end{align*}
\end{proposition}
\begin{proof}
See Appendix~\ref{app:Prop4}.
\end{proof}

We now present a composite relation for the four quantities defined in \eqref{eq-quants}, by defining the vector $V_k$ as follows:
\begin{align}\label{eq-vk}
\!\!V_k=\Big(\|\hat x_k -x^*\|,D(\bx_k,\!\phi_k),S(\by_k,\!\pi_k),\|\bx_{k}\!-\!\bx_{k-1}\|\Big)^{\!\T}\!\!\!.
\end{align}
For the vector $V_k$ we  have the following result.
			
\begin{proposition}\label{prop-main}
Let Assumption~\ref{asm-graphs}, Assumption~\ref{asm-functions}, and Assumption~\ref{asm-strconv} hold. Then, for the iterates produced by the accelerated $AB$/Push-Pull method in~Algorithm 1 with the step-size $\alpha\in(0,2(nL)^{-1})$, we have
\[V_{k+1}\le M_k(\a,\b,\g)V_k\qquad \hbox{for all $k\ge0$},\]
where the $ij$-th element $m_k^{i,j}$ of the matrix $M_k(\a,\b,\g)$ are given as follows:
\begin{align*}
&m_k^{1,1}=q_k(\a),~m_k^{1,2}=\a L\sqrt{n}\varphi_k,~m_k^{1,4}=\b\!+\!\g(1\!+\!\a L\sqrt{n}),\\
&m_k^{1,3}=\a,~m_k^{2,1}=\a L\sqrt{n}\varphi_k,~m_k^{2,2}=c_k+\a L\sqrt{n}\varphi_k,\\
&m_k^{2,3}=\a,~~m_k^{2,4}=\b+\g (c_k+\a L\sqrt{n}),\\
&m_k^{3,1}=Lr_k(1+\g)m_k^{4,1},~~m_k^{3,2}=Lr_k(1+\g)m_k^{4,2},\\
&m_k^{3,3}=\tau_k\!+\!Lr_k(1\!+\!\g)m_k^{4,3},~m_k^{3,4}=L r_k\g \!+\! Lr_k(1\!+\!\g)m_k^{4,4},\\
&m_k^{4,1}=\a L\sqrt{n}\varphi_k,~~m_k^{4,2}=c_k\varphi_{k+1}+\varphi_k+\a L\sqrt{n}\varphi_k,\\
&m_k^{4,3}=\a,~~m_k^{4,4}=\b+\g\sqrt{n}(1+\a L).
\end{align*}
\end{proposition}
\begin{proof}
The stated relation follows directly from Proposition~\ref{prop-waverx}, Proposition~\ref{prop-xcontract}, Proposition~\ref{prop-xdiff}, and Proposition~\ref{prop-ycontract}.
\end{proof}

\begin{remark}
From Proposition~\ref{prop-main}, to prove that $V_k$ tends to 0 at a geometric rate, it is sufficient to show that 
\[M_k(\a,\b,\g)\le M(\a,\b,\g),\]
for some matrix $M(\a,\b,\g)$, where the preceding inequality is to be understood entry-wise. Then, we select an appropriate step-size $\alpha$ in the range $(0,2(nL)^{-1})$ and the acceleration parameters $\b$, $\g$, such that the eigenvalues of $M(\a,\b,\g)$ lie inside the unit circle, i.e., the spectral radius of $M(\a,\b,\g)$ is less than $1$.
\end{remark}

We now determine an upper bound matrix $M(\alpha, \beta, \gamma)$ for $M_k(\alpha, \beta, \gamma)$. To do this, we define upper bounds for the constants $c_k, \tau_k, r_k,$ and $\varphi_k$ defined in \eqref{eq-const-all} as $c \in (0,1)$, $\tau \in (0,1)$, $r$, and $\varphi$, respectively, i.e.,
\begin{align}\label{eq-maxconst}
	\!\!\!\max_{k\ge0}c_k\le c,~ \max_{k\ge0}\tau_k\le \tau,~ \max_{k\ge0}r_k\le r,~ \max_{k\ge0} \varphi_k \le \varphi.
\end{align}
For the quantity $q_k(\a)$ in Lemma~\ref{prop-waverx}, when $\a\in(0,2(nL+n\mu)^{-1})$, we have $q_k(\a)=1-\a n \min(\pi_k)\mu <1$. Let $\sigma$ be a lower bound for $\min(\pi_k)$, $k\ge0$, i.e., $\sigma \le \min_{k\ge0}\{\min(\pi_k)\}$. Note that Lemma~\ref{lem-bmatrices} provides such a lower bound applicable to any sequence of strongly connected graphs $\{\mathbb{G}_k\}$. Better lower bounds can be obtained when the graphs have more specific structures. We have the following upper bound for $q_k(\a)$:
\begin{align}\label{eq-maxconst2}
	\max_{k\ge0}q_k(\a) \le 1-\a n \sigma \mu\in (0,1).
\end{align}
Using the upper-bounds given in~\eqref{eq-maxconst} and~\eqref{eq-maxconst2}, for $\a\in(0,2(nL+n\mu)^{-1})$, we have $M_k(\a,\b,\g)\le M(\a,\b,\g)$, for all $k \ge 0$, with the matrix $M(\a,\b,\g)$ given by
\begin{align}\label{eq-gmatrixm}
	\!\!\!\!\!\left[\!\begin{array}{cccc}
	 \!\!1\!-\!\a n \sigma \mu\!\! & \a L \sqrt{n} \varphi & \a & \b\!+\!\g(1\!+\!\a L\sqrt{n})\!\! \cr
	 \hbox{}\cr
	 \a L\sqrt{n}\varphi & c + \a L\sqrt{n}\varphi  & \a & \!\!\b \!+\!\g (c\!+\!\a L\sqrt{n})\!\!\cr
	 \hbox{}\cr
	 u_1 & u_2 & \!\tau \!+\! u_3\!& Lr\g\!+\!u_4\cr
	 \hbox{}\cr
	 \!\a L\sqrt{n}\varphi & \!\!\!(c\!+\!1)\varphi \!+\! \a L\sqrt{n}\varphi\!\!\!  & \a & \!\!\b\!+\!\g\sqrt{n}(1\!+\!\a L)\!
	 \end{array}\!\right]\!\!\!\!\!
\end{align}
where the third row of the matrix $M(\a,\b,\g)$ is co-linear with the fourth row, i.e.,
\[(u_1,u_2,u_3,u_4)=Lr(1+\g)[M(\a,\b,\g)]_{4,:}.\]
We next define constants $\eta_i$, $i\in [6]$, as follows:
\begin{subequations}\label{eq-const_eta}
\begin{align}
	&\eta_1 = (1-\tau)(1-c)n\sigma\mu\label{eq-const_eta1},\\
	&\eta_2 = (1-\tau)(n\sigma\mu L\sqrt{n}\varphi+L^2n\varphi^2),\\
	&\eta_3 = Lr[(1+c)\varphi+1-c](n\sigma\mu+L\sqrt{n}\varphi),\\
	&\eta_4 = (1-\tau)[(1+c)\varphi+1-c]\label{eq-const_eta4},\\
	&\eta_5 = \eta_1(\sqrt{n}\!-\!c)\!+\![n\sigma\mu(1\!+\!c\!+\!L\sqrt{n})\!+\!2L\sqrt{n}\varphi]\eta_4\label{eq-const_eta5},\\
	&\eta_6 = (1+\g c-\g\sqrt{n})\eta_2
	+(1+\g)\eta_3 +L^2n\varphi\eta_4.\label{eq-const_eta6}
\end{align}
\end{subequations}
We now present the main result of this paper, which states that the accelerated $AB$/Push-Pull algorithm (Algorithm 1) converges to the global minimizer at a linear rate.
\begin{theorem}\label{theo-main}
\rev{Let Assumption~\ref{asm-graphs}, Assumption~\ref{asm-functions}, and Assumption~\ref{asm-strconv} hold.} Consider the iterates produced by Algorithm 1, the notations in \eqref{eq-maxconst}-\eqref{eq-maxconst2} and the constants $\eta_i,~i\in [6]$, defined in \eqref{eq-const_eta}. If the step-size $\a>0$ and the acceleration parameters 
$\b\ge0$ and $\g\ge 0$ are chosen such that 
\begin{align}\label{eq:alpha-range}
\begin{cases} 
	\a \le \min \Bigg\{\dfrac{1-c}{L\sqrt{n}\varphi},\dfrac{1-\tau}{Lr},\dfrac{\eta_1-\kappa\eta_5}{\eta_6},\dfrac{2}{n(L+\mu)}\Bigg\},\!\!\!\!\!\!\!\\
	\max\{\b,\g\}<\dfrac{\eta_1}{\eta_5},\\
	\b + \g \sqrt{n}< 1,
\end{cases} 
\end{align}
then $\rho_{M}<1$ where $\rho_{M}$ is the spectral radius of $M(\a,\b,\g)$  and, thus, $\lim_{k\to\infty} \|x_k^i -x^*\|=0$ with a linear convergence rate of the order of $\mathcal{O}\Big(\rho_M^k\Big)$ for all $i\in [n]$.
\end{theorem}

\begin{proof}
Recall that by Proposition~\ref{prop-main}, we have \[V_{k+1}\le M_k(\a,\b,\g)V_k,\qquad \text{for all}~ k\ge 0.\] 
With the matrix $M(\a,\b,\g)$ defined as in \eqref{eq-gmatrixm}, we obtain
\begin{equation}\label{eq-vkrel}
	V_{k+1}\le M(\a,\b,\g)V_k, \qquad \hbox{for all $k\ge0$}.
\end{equation}
Thus, $\|\hat x_k -x^*\|,~D(\bx_k,\phi_k),~\|\bx_{k}\!-\!\bx_{k-1}\|$, and $S(\by_k,\pi_k)$ all converge to $0$ linearly at rate $\mathcal{O}\Big(\rho_M^k\Big)$ if the spectral radius $\rho_{M}$ of $M(\a,\b,\g)$ satisfies $\rho_{M}<1$. By Lemma~8 of~\cite{pshi21}, we will have $\rho_{M}<1$ if all diagonal entries of $M(\a,\b,\g)$ are less than 1 and $\det(\mbi-M(\a,\b,\g))>0$ where
\begin{align*}
	&\det(\mbi-M(\a,\b,\g)) \\
	= &\a\eta_1(1+\g c-\g\sqrt{n})-\a^2(1+\g c-\g\sqrt{n})\eta_2 -\a^2 (1+\g)\eta_3\\
	-&\a[n\sigma\mu(\b+\g c+\g L\sqrt{n})+L\sqrt{n}\varphi(\b+\g)]\eta_4-\a^2 L^2n\varphi\eta_4,
\end{align*}
with positive constants $\eta_1>0$, $\eta_2>0$, $\eta_3>0$ and $\eta_4>0$ defined as in \eqref{eq-const_eta1}--\eqref{eq-const_eta4}.
Let $\kappa = \max\{\b,\g\}$, we can further simplify the determinant as follows
\begin{align*}
	\det(\mbi-M(\a,\b,\g)) = \a\left[(\eta_1-\kappa\eta_5)-\a\eta_6\right],
\end{align*}
where $\eta_5$ and $\eta_6$ are positive constants defined in equations (\ref{eq-const_eta5}) and (\ref{eq-const_eta6}), given that $\g \sqrt{n} < 1$. 
Hence, we need to choose $\a\in(0,2(nL+n\mu)^{-1})$ and $\b\ge0$, $\g\ge0$, so that the following conditions are satisfied
    \[\begin{cases} 
    c+\a L\sqrt{n}\varphi<1, \quad \tau + \a Lr< 1,\\
    \b + \g \sqrt{n}< 1,\quad (\eta_1-\kappa\eta_5)-\a\eta_6>0.
    \end{cases}\]
which yields the range in \eqref{eq:alpha-range}.
    \end{proof}

We note that Theorem~\ref{theo-main} includes the case when $\b=\g=0$, thus recovering the convergence rate of the $AB$/Push-Pull for time-varying graphs without accelaration, which has been shown in~\cite{Saadatniaki2020} and, recently, in~\cite{Angelia2022AB} (with explicit bounds on the stepsize selection). Theorem~\ref{theo-main} also encompasses the results of acceleration for static directed graphs where the graph remains unchanged over time, $\bbG_k=\bbG$, for all time steps, $k>0$. This includes the case when $\g=0$ from \cite{Xin2020Heavy}, the missing convergence analysis for $\b=0$ in \cite{Xin2019Nes}, and the case $\b=\g>0$ presented in \cite{Huaqing2021}.

\begin{remark}
The convergence analysis for $C$-strongly-connected graph sequences is performed similarly to our analysis above by utilizing the results in Remark~\ref{rem-CgraphsAB}, and by recognizing that contractions resulting from row- and column-stochastic matrices occur after $k=C$.
\end{remark}

\section{Numerical Simulations}\label{sec:simulation}
In this section, we evaluate the performance of the proposed algorithm by testing it on several real-world datasets, and assessing its accuracy and efficiency.
We compare the results between the $AB$/Push-Pull algorithm (ABPP) and its variants using different acceleration techniques, including the heavy-ball momentum (ABPP-m), Nesterov momentum (ABPP-N), and the combination of the two momentum techniques (ABPP-mN). This comparison is sufficient as the performance of the $AB$/Push-Pull algorithm and other existing algorithms, such as Push-DIGing \cite{nedic2017achieving} and subgradient-push \cite{nedic2015distributed}, \rev{as discussed in our introduction}, have already been evaluated (see, for example,  \cite{Angelia2022AB,Saadatniaki2020}). 

In our simulations, all the communication graphs are directed, time-varying, and have self-loops. To ensure the graphs are strongly connected, a directed cycle linking all agents is established at each iteration. \rev{Utilizing time-varying directed communication graphs proves to be highly practical in numerous scenarios characterized by dynamic communication networks among agents, where the flow of information or commands between agents can be directed. These scenarios often arise due to various factors such as communication delays, user mobility, and the influence of straggler effects.}


\subsection{Distributed Ridge Regression}\label{sec:simulation-Least-squares}
Consider a sensor fusion problem, as described in \cite{Jinming2018,pshi21}. The goal of the sensor fusion problem is to estimate an unknown parameter $x$ by utilizing data from $n$ sensors. Each sensor $i$ has a measurement matrix $H_i \in \mathbb{R}^{s\times p}$, and a noisy observation $z_i=H_ix+\omega_i \in \mathbb{R}^{s}$ of $x$, where $\omega_i$ represents the noise. The resulting sensor fusion problem is given by the following minimization problem:
\[\min_{x\in \mathbb{R}^p}~ \sum\limits_{i=1}^n \Big(\|z_i-H_ix\|^2+\lambda_i \|x\|^2\Big),\]
where  $\lambda_i>0$ is the regularization parameter for the local cost function of sensor $i$.

\begin{figure}[h]
    \vspace{-0.2cm}
	\centering
	\hspace*{-1.1em}
		\subfigure[Sensor fusion problem]{
	 \includegraphics[width=0.247\textwidth]{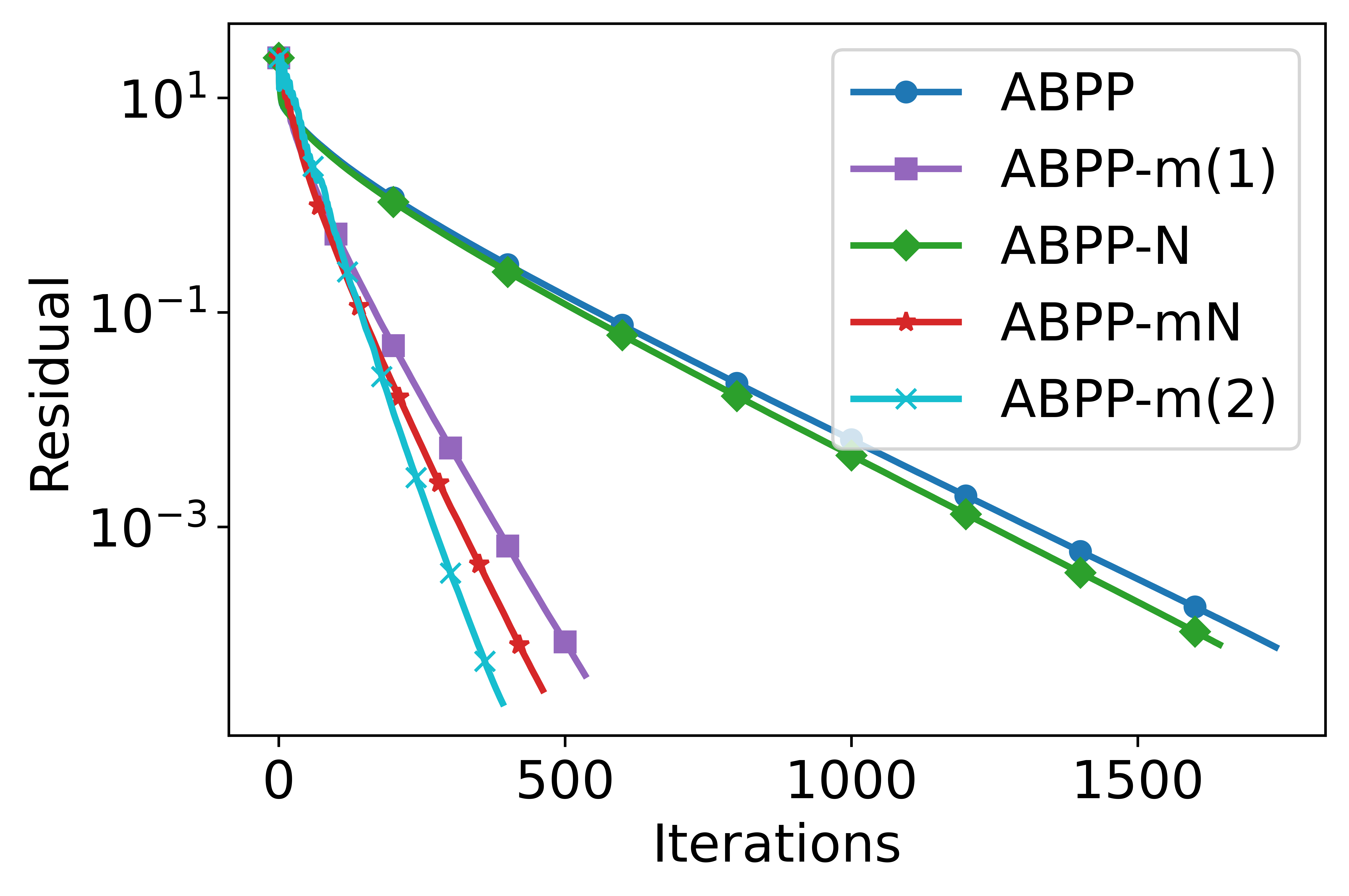}
	    \label{fig:Sensor}  } 
	   \hspace*{-1.5em} 
		 \subfigure[Diabetes dataset]{
	     \includegraphics[width=0.247\textwidth]{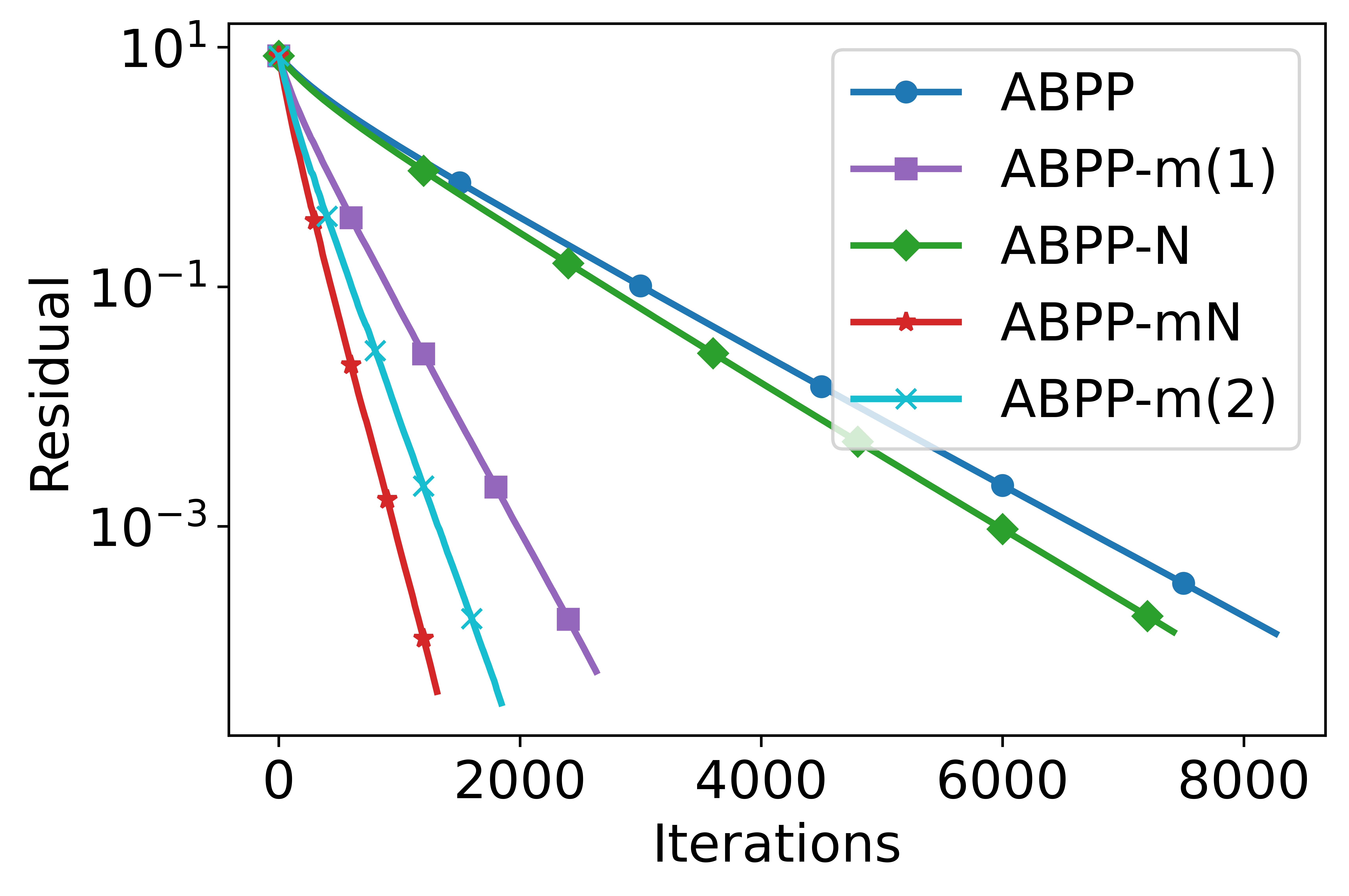}
	   \label{fig:Diabetes}
	} \hspace*{-1.8em}
 \vspace{-0.2cm}
	\caption{Residual plots}\label{fig:residual}
 \vspace{-0.3cm}
\end{figure}

We follow the setup in \cite{pshi21} where $n=20$, $p=20$, and $s=1$ are chosen to make the local cost function ill-conditioned, requiring coordination among agents for fast convergence. The measurement matrix $H_i$ is generated from a uniform distribution in the unit $\re^{s\times p}$ space and, then, normalized so that its Lipschitz constant is equal to $1$. The noise $\omega_i$ follows an i.i.d.\ Gaussian process with zero mean and unit variance $\mathcal{N}(0,1)$. The regularization parameter is chosen as $\lambda_i=0.01$ for all $i\in [n]$ to ensure the strong convexity of the loss function. Figure~\ref{fig:Sensor} illustrates the accelerated linear convergence of the proposed algorithm using different acceleration techniques \rev{(ABPP, ABPP-m(1), ABPP-N, ABPP-mN), using $\a=0.25$, $\b=0.7$ and $\g=0.05$. When Nesterov momentum is not used, a larger value of $\b=0.8$ may be selected (ABPP-m(2)).} The plot measures performance based on the residual between the iterates $x_k^i$, $i\in[n]$, at time step $k$ and the optimal value $x^*$, given by $\tfrac{1}{n}\sum_{i=1}^n\|x_k^i-x^*\|$.

\begin{figure}[h]
\vspace{-0.1cm}
	\centering
	\hspace*{-1.7em}
	\includegraphics[width=0.24\textwidth]{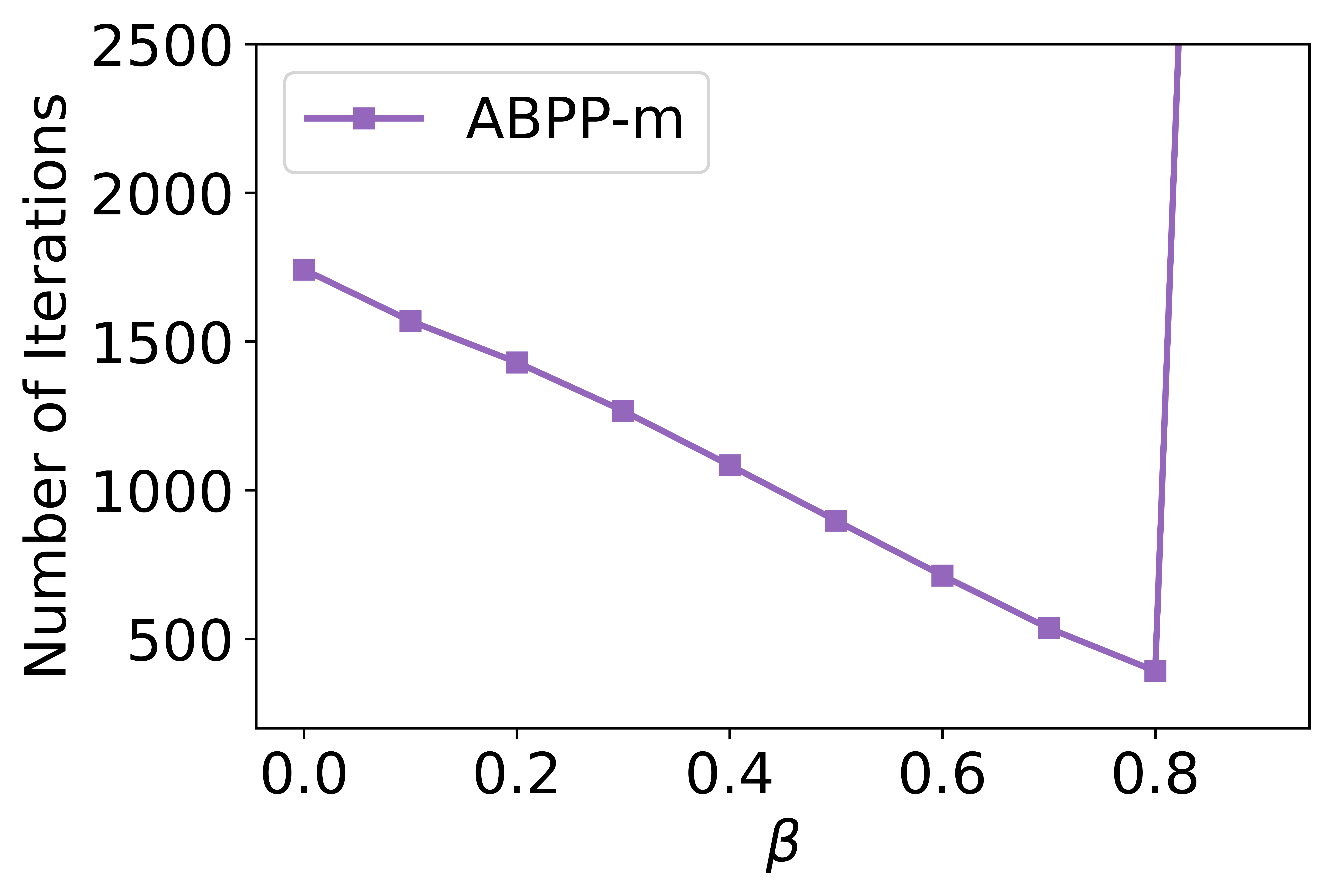}
	\includegraphics[width=0.24\textwidth]{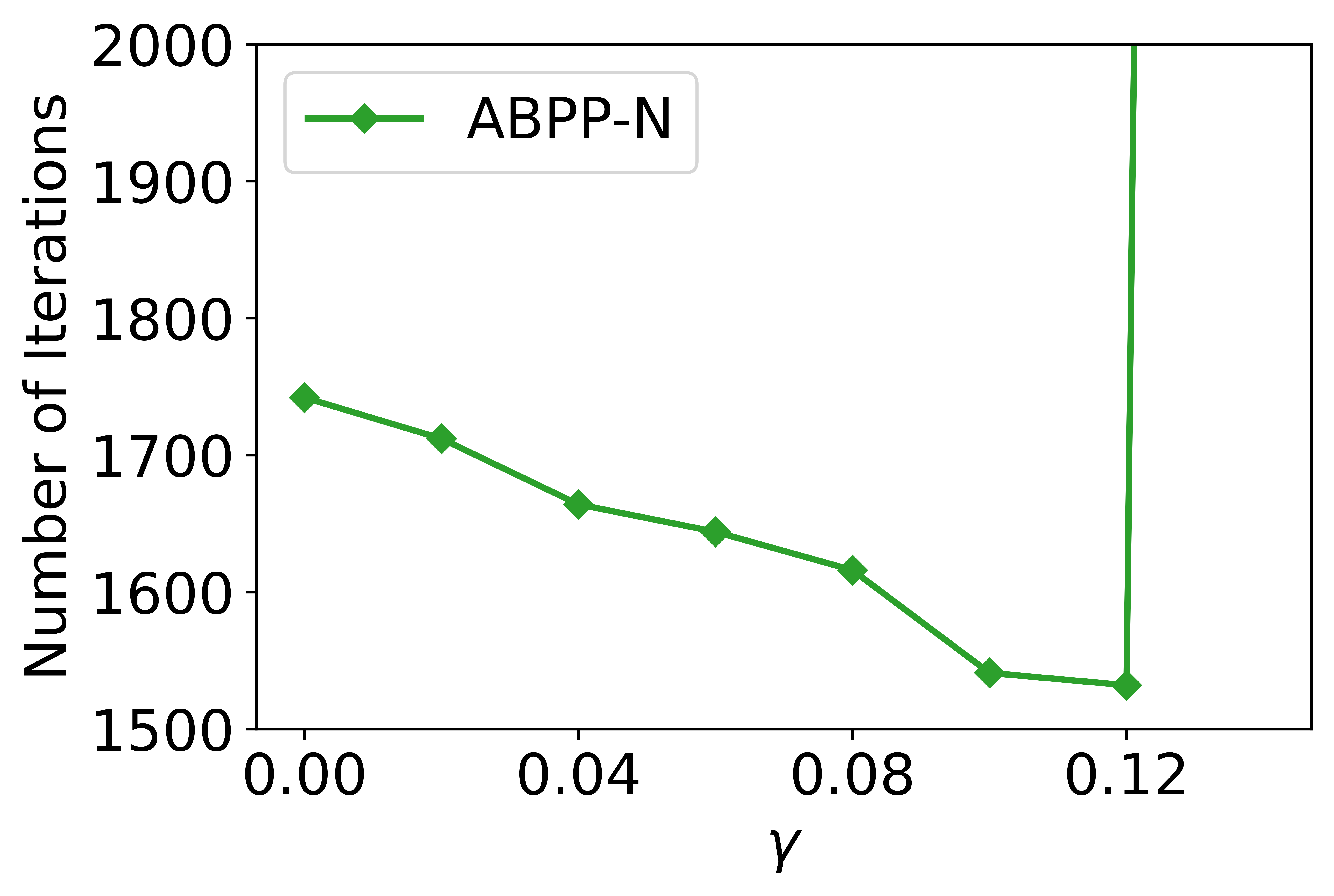}
	\hspace*{-1.7em}
 \vspace{-0.2cm}
	\caption{Effects of varying momentum parameters ($\a=0.25$)}\label{fig:effects}
 \vspace{-0.5cm}
\end{figure}

We also test how the changes in the momentum parameters $\b$ and $\g$ affect the convergence rate, as shown in Figure~\ref{fig:effects}. The results imply that the algorithm converges faster with higher momentum parameter values. Also, the range of parameters satisfies the conditions in~\eqref{eq:alpha-range}. Note that the value of $\g$ is much smaller due to the condition $\b+\g\sqrt{n}<1$ in~\eqref{eq:alpha-range}.


\subsection{Distributed Logistic Regression (L2-Regularization)}\label{sec:simulation-Logistic}
In this experiment,  we examine binary classification problems using real-world datasets to evaluate the performance among the different acceleration technique over a time-varying directed network. We consider a total of $N$ labeled data point for training, with each node $i$ possessing a local batch of $m_i$ training samples. The $j$-th sample at node $i$ is a tuple $\{b_{ij},y_{ij}\} \!\!\!\subseteq \!\!\re^p \!\times \!\!\{
+1,\!-1\}$. To construct an estimate \rev{$x\!\!=\!\![x_0,x_{1:}^\T]^\T\!\!\in\! \re^{p+1}$ of the coefficients, where $x_{1:}\!\!=\!\![x_1,\ldots,x_p]^\T\!$}, we will use the principle of maximum likelihood and define the local logistic regression cost function $f_i$ at node $i$ as:
\[\rev{f_i(x) = \dfrac{1}{m_i}\sum_{j=1}^{m_i}\ln \Big[1+\exp\Big\{-(x_{1:}^\T b_{ij}+x_0)y_{ij}\Big\} \Big] +\dfrac{\lambda}{2}\|x\|^2\!\!,}\]
which is smooth and strongly convex due to the inclusion of the L2-regularization. The nodes cooperate to solve the following decentralized optimization problem:
\[\rev{\min_{x\in \mathbb{R}^{p+1}}~ \frac{1}{n}\sum\limits_{i=1}^n f_i(x).}\]
We examine the performance of the proposed algorithm using two datasets.

\textbf{Pima Indians Diabetes Dataset:}
We evaluate the performance of our algorithm using the Diabetes dataset, which consists of $N=700$ training samples and $68$ test samples. The dataset is divided among $n=7$ agents, with each agent having $m_i=100$ samples. A regularization parameter of $\lambda=0.001$ is used and the algorithm stops when the consensus error among agents is less than $10^{-7}$. The accuracy on the test set is $79.41\%$ for all $4$ algorithms. Figure~\ref{fig:Diabetes} illustrates the accelerated linear convergence of the proposed algorithm using different acceleration techniques for $\a=0.5$, \rev{$\b=0.7$ (ABPP-m(1), ABPP-mN)} and $\g=0.1$ \rev{while $\b=0.8$ for ABPP-m(2)}.

\begin{figure}[ht!]
    \centering
    \includegraphics[width=0.35\textwidth]{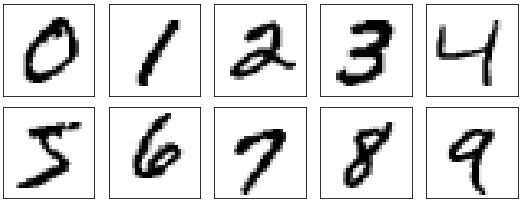}
    \caption{Samples from MNIST Dataset}\label{fig:MNIST_sample}
    \vspace{-0.3cm}
\end{figure}

\textbf{MNIST Dataset:}
The task at hand is to perform digit classification on the MNIST dataset. Figure~\ref{fig:MNIST_sample} shows a part of randomly selected samples, where each image is featured as a $784$-dimensional vector. We use $2000$ training samples and $1000$ test samples. The problem is divided among $10$ agents, with each agent handling $200$ samples. The regularization parameter is set to $\lambda=0.001$, and the algorithm terminates when the consensus error among agents is below $10^{-3}$.

\begin{figure}[ht!]
    \centering
    \includegraphics[width=0.2\textwidth]{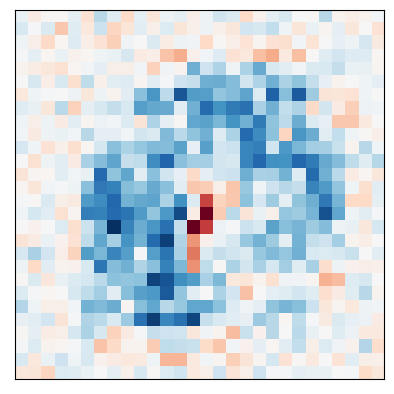}
    \caption{Heatmap of Coefficients for Logistic Regression in Classifying Zero and Non-Zero Digits using ABPP-mN}\label{fig:Coefficients}
\end{figure}

As a sanity check of accuracy, we visually examine the coefficients generated by the proposed algorithm for the binary classification task of hand-written digit 0 versus non-zero digits. The visualization, shown in Figure~\ref{fig:Coefficients}, highlights the most important features identified by the algorithm for classifying a digit as $0$. The blue pixels indicate the highest probability of a digit being classified as $0$, while the red pixels indicate the lowest probability. As seen in Figure~\ref{fig:Coefficients}, the blue pixels form the shape of a $0$, with more red pixels in the center, indicating that these pixels are less likely to be shaded in images of a $0$.

We now assess the performance of the four algorithms in classifying hand-written digits. The tasks are to classify hand-written digits $\{1, 2\}$ and classify  hand-written digits $\{3, 5\}$. The results of the numerical experiments are shown in Table~\ref{tab:Accuracy} which include the number of iterations, \rev{the computational time} and the accuracy on the test set. 

\begin{table}[h!]
\centering
\begin{tabular}{|l|ccc|ccc|}
\hline
\multirow{2}{*}{}            & \multicolumn{3}{c|}{Classify $\{1,2\}$}        & \multicolumn{3}{c|}{Classify $\{3,5\}$}        \\ \cline{2-7} 
                             & \multicolumn{1}{c|}{\!\!\!\#\! Iterations\!\!\!\!} & \multicolumn{1}{c|}{\!\!\!\rev{Time (s)}\!\!\!\!}  & \!\!\!\!Accuracy\!\!\!\!  & \multicolumn{1}{c|}{\!\!\!\#\! Iterations\!\!\!\!} & \multicolumn{1}{c|}{\!\!\!\!\rev{Time (s)}\!\!\!\!}  & \!\!\!\!Accuracy\!\!\!\!  \\ \hline
\!\!\!\!ABPP& \multicolumn{1}{c|}{$1326$}& \multicolumn{1}{c|}{$75.5$}        & $98.3\%$ & \multicolumn{1}{c|}{$3287$}& \multicolumn{1}{c|}{$246.0$}        & \!\!$95.49\%$\!\! \\ \hline
\!\!\!\!ABPP--m & \multicolumn{1}{c|}{$1008$}& \multicolumn{1}{c|}{$51.9$}        & $98.6\%$ & \multicolumn{1}{c|}{$1704$}& \multicolumn{1}{c|}{$137.4$}        & \!\!$95.30\%$\!\! \\ \hline
\!\!\!\!ABPP--N    & \multicolumn{1}{c|}{$1380$}& \multicolumn{1}{c|}{$95.8$}        & $98.3\%$ & \multicolumn{1}{c|}{$3261$}& \multicolumn{1}{c|}{$266.5$}        & \!\!$95.49\%$\!\! \\ \hline
\!\!\!\!ABPP--mN\!\!\!\! & \multicolumn{1}{c|}{$944$}& \multicolumn{1}{c|}{$47.6$}         & $98.6\%$ & \multicolumn{1}{c|}{$1606$}& \multicolumn{1}{c|}{$117.8$}        & \!\!$95.70\%$\!\! \\ \hline
\end{tabular}
\caption{Performance for $\alpha=0.01$, $\beta=0.3$ and $\gamma=0.01$. }\label{tab:Accuracy}
\vspace{-0.5cm}
\end{table}


Overall, the results demonstrate that the proposed algorithm with acceleration techniques substantially enhances the convergence rate while having comparable performance to the $AB$/Push-Pull algorithm. \rev{The heavy-ball momentum, in particular,  significantly improves the convergence, and incorporating both heavy-ball and Nesterov momentum may be beneficial in some cases (Figure~\ref{fig:residual})}.
The Nesterov parameter $\g$ is influenced by the number of agents in the network, as it is multiplied by $\sqrt{n}$ (see \eqref{eq:alpha-range}). The heavy-ball parameter $\b$, on the other hand, can be set to a larger value for faster convergence. By considering different values for the Nesterov and heavy-ball acceleration parameters $\b$ and $\g$, our proposed algorithm offers greater flexibility and the potential for faster convergence.

\section{Conclusion}\label{sec:conclusion}
In this paper, we propose a novel approach for solving distributed optimization problems over time-varying directed networks. The proposed algorithm incorporates acceleration techniques to improve the performance of the $AB$/Push-Pull algorithm. Theoretical analysis is provided to prove linear convergence to the optimal solution under certain conditions. Additionally, explicit bounds for the step-size and momentum parameters are derived based on the properties of the cost functions and network structure. The numerical results demonstrate the benefits of the proposed acceleration techniques on the $AB$/Push-Pull algorithm. An interesting open direction is to theoretically analyze the acceleration over the $AB$/Push-Pull algorithm.

\bibliographystyle{IEEEtran}
\bibliography{ref_short}

\vspace{-0.5cm}
\begin{IEEEbiography}[{\includegraphics[width=1in,height=1.1in,clip,keepaspectratio]{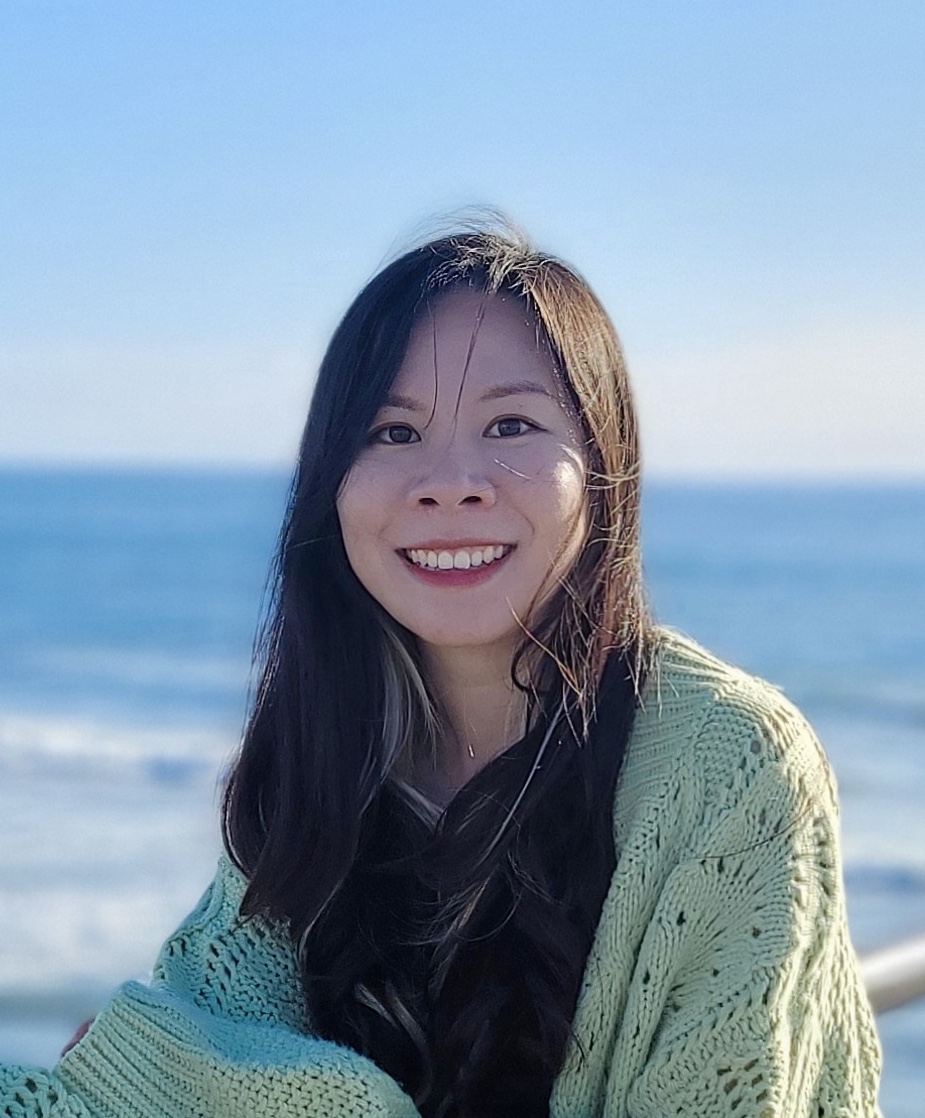}}]
	{Duong Thuy Anh Nguyen} (Graduate Student Member, IEEE) received the M.Sc. degree in Applied Mathematics from the University of Louisiana at Lafayette, LA, USA in 2019. She is currently pursuing the Ph.D. degree with the School of Electrical, Computer and Energy Engineering at Arizona State University, AZ, USA. Her current research interests include distributed optimization, operations research and game theory. The research focuses on developing mathematical models for decision-making under uncertainty, fair and privacy-preserving mechanism designs and distributed algorithms for large-scale network in multi-agent systems. 
	Research applications include cloud/edge computing, electric vehicles, non-cooperative games over communication networks.
	\end{IEEEbiography}
\vspace{-0.5cm}	
\begin{IEEEbiography}[{\includegraphics[width=1in,height=1.25in,clip,keepaspectratio]{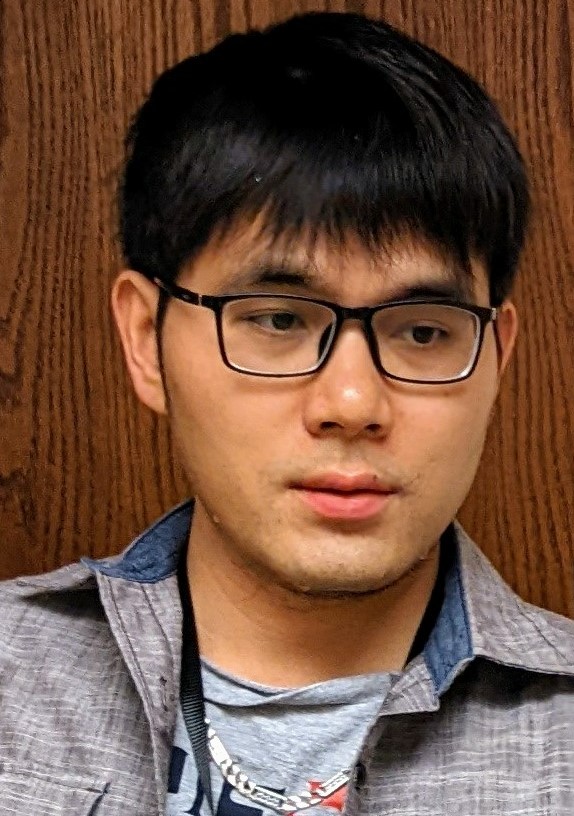}}]
    {Duong Tung Nguyen} received the Ph.D. degree in electrical and computer engineering from the University of British Columbia, BC, Canada. He is currently an assistant professor in the School of Electrical, Computer and Energy Engineering at Arizona State University, AZ, USA. His research lies at the intersection of operations research, AI, economics, and engineering, with a focus on developing new mathematical models and techniques for decision-making and economic analysis of large-scale networked systems such as cloud/edge computing, smart grids, intelligent transportation, and crowdsourcing. 
	\end{IEEEbiography}
\vspace{-0.5cm}
\begin{IEEEbiography}[{\includegraphics[width=1in,height=1.25in,clip,keepaspectratio]{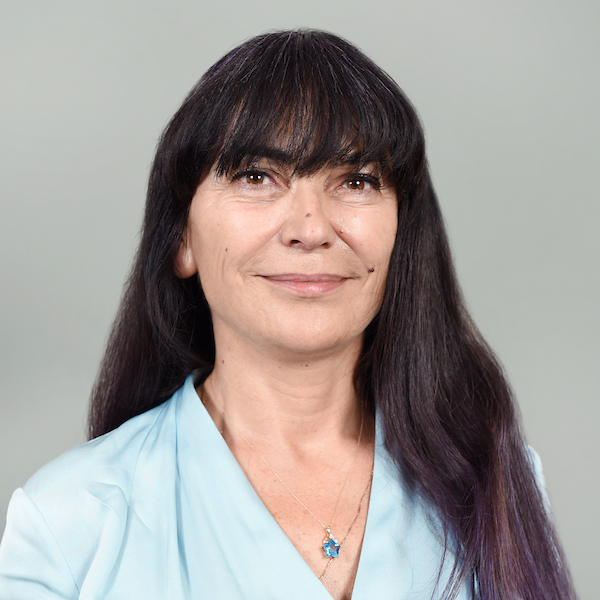}}]
	{Angelia Nedi\'c} has a Ph.D.\ from Moscow State University, Moscow, Russia, in Computational Mathematics and Mathematical Physics (1994), and a Ph.D.\ from Massachusetts Institute of Technology, Cambridge, USA, in Electrical and Computer Science Engineering (2002). She has worked as a senior engineer in BAE Systems North America, Advanced Information Technology Division at Burlington, MA. Currently, she is a faculty member of the School of Electrical, Computer and Energy Engineering at Arizona State University at Tempe. Prior to joining Arizona State University, she has been a Willard Scholar faculty member at the University of Illinois, Urbana-Champaign. 
    She is a recipient (jointly with her co-authors) of the Best Paper Award at the Winter Simulation Conference 2013 and the Best Paper Award at the International Symposium on Modeling and Optimization in Mobile, Ad Hoc and Wireless Networks (WiOpt) 2015.  Her general research interest is in optimization, large scale complex systems dynamics, variational inequalities, and games. 
	\end{IEEEbiography}

\appendix
\subsection{Proof of Proposition~\ref{prop-waverx}}\label{app:Prop1}
\begin{proof}
    Under Assumption~\ref{asm-strconv}, the unique minimizer $x^*$ of $f(x)$ over $x\in\re^p$ exists. 
	By subtracting $x^*$ from both sides of the relation for $\hat x_k$ in Lemma~\ref{lem-weightedavg}, and by
	adding and subtracting $\sum_{i=1}^n [\phi_{k+1}]_i \a n [\pi_k]_i \nabla f(\hat x_k)$, we obtain:
	\begin{align*}
	\hat x_{k+1}-x^*&=\hat x_k - x^* - \sum_{i=1}^n [\phi_{k+1}]_i\a n [\pi_k]_i\nabla f(\hat x_k) \\
	&+\a \sum_{i=1}^n [\phi_{k+1}]_i\Big(n [\pi_k]_i \nabla f(\hat x_k) -  y_k^i\Big)\\
	&+\sum_{i=1}^n (\b[\phi_{k+1}]_i+\g[\phi_{k}]_i)(x_k^i-x_{k-1}^i).
	\end{align*}
	As a result of the convexity of the norm and the stochastic nature of $\phi_{k+1}$, we can deduce that:
	\begin{align}\label{eq-proofweightx}
	\|\hat x_{k+1}-x^*\|&\le \sum_{i=1}^n [\phi_{k+1}]_i\|\hat x_k - x^* - \a n [\pi_k]_i\nabla f(\hat x_k) \| \nonumber\\
	&+\a \sum_{i=1}^n[\phi_{k+1}]_i \|y_k^i-n [\pi_k]_i \nabla f(\hat x_k)\| \nonumber\\
	&+\sum_{i=1}^n (\b[\phi_{k+1}]_i+\g[\phi_{k}]_i)\|x_k^i-x_{k-1}^i\|.
	\end{align}
	For a step-size $\a$ satisfying $\a\in  (0,\tfrac{2}{n[\pi_k]_iL})$, for all $i\in [n]$, by Lemma~\ref{lem-contraction}, it follows that 
	\[\|\hat x_k - x^* - \a n [\pi_k]_i\nabla f(\hat x_k) \| \le q_{i,k}(\a) \|\hat x_k - x^*\|,\]
	with $q_{i,k}(\a)=\max\{|1-\a n [\pi_k]_i\mu|,|1-\a n [\pi_k]_i L|\}$.
	
 We have the following estimate for the first term on the right-hand-side of \eqref{eq-proofweightx}:
	\begin{align*}
    &\sum_{i=1}^n[\phi_{k+1}]_i\|\hat x_k - x^*-\a n [\pi_k]_i\nabla f(\hat x_k)\| \\
    \le &\sum_{i=1}^n[\phi_{k+1}]_i q_{i,k}(\a) \|\hat x_k - x^*\| \le q_k(\a) \|\hat x_k - x^*\|,
    \end{align*}
    using the stochasticity of $\phi_{k+1}$ and $q_k(\a) = \max\{|1-\a n \min(\pi_k)\mu|,|1-\a n \min(\pi_k)L|\}$.
	Since $\max(\phi_{k+1})\le 1$, to estimate the second term on the right-hand-side in~\eqref{eq-proofweightx}, we factor out $[\pi_k]_i$ (which is positive by Assumption~\ref{asm-graphs} and Lemma~\ref{lem-bmatrices}), as follows
     \begin{align*}
    &~~~\sum_{i=1}^n[\phi_{k+1}]_i\|y_k^i\!-\!n [\pi_k]_i \nabla f(\hat x_k)\| \\
    &\le  \!\sum_{i=1}^n\|y_k^i\!-\!n [\pi_k]_i \nabla f(\hat x_k)\| = \!\sum_{i=1}^n[\pi_k]_i \Bigg\|\dfrac{y_k^i}{[\pi_k]_i}\!-\!n \nabla f(\hat x_k)\Bigg\|\\
    &\le \sum_{i=1}^n[\pi_k]_i \Bigg\|\dfrac{y_k^i}{[\pi_k]_i}-\sum_{\ell=1}^n y_k^\ell\Bigg\| +\sum_{i=1}^n[\pi_k]_i\Bigg\|\sum_{\ell=1}^n y_k^\ell-n\nabla f(\hat x_k)\Bigg\|\cr
    &\!\le\! \sqrt{\sum_{i=1}^n[\pi_k]_i \Bigg\|\dfrac{y_k^i}{[\pi_k]_i}\!-\!\sum_{\ell=1}^n y_k^\ell\Bigg\|^2} \!+\Bigg\|\sum_{\ell=1}^n y_k^\ell\!-\!n\nabla f(\hat x_k)\Bigg\|
    \\
    &\le S(\by_k,\pi_k) \!+\Bigg\|\sum_{\ell=1}^n y_k^\ell\!-\!n\nabla f(\hat x_k)\Bigg\|,
    \end{align*}
    where we add and subtract $\sum_{\ell=1}^n y_k^\ell$, and use the triangle inequality to obtain the second inequality. We use the fact that the vector sequence $\{\pi_k\}$ is stochastic to derive the third inequality (see Lemma~\ref{lem-amatrices}), and the last inequality follows from the definition of the $S$-quantity in \eqref{eq-x-S-quants}. 
    We now estimate the last term in the preceding relation. By Lemma~\ref{lem-sumgrad} we have
	$\sum_{\ell=1}^n y_k^\ell=\sum_{\ell=1}^n \nabla f_\ell(s_k^\ell)$;
	hence, in view of $\nabla f=\frac{1}{n}\sum_{\ell=1}^n \nabla f_\ell$, it follows that
	\begin{align*}
	&\Bigg\|\sum_{\ell=1}^n y_k^\ell-n\nabla f(\hat x_k)\Bigg\|
	= \Bigg\|\sum_{\ell=1}^n \Big(\nabla f_\ell(s_k^\ell) - \nabla f_\ell(\hat x_k)\Big)\Bigg\|\\
	&\le \sum_{\ell=1}^n \|\nabla f_\ell(s_k^\ell) -  \nabla f_\ell(\hat x_k)\|\le L\sum_{\ell=1}^n \|s_k^\ell -  \hat x_k\| \\
 &= L\sqrt{n}\varphi_k D(\bx_k,\phi_k)+L\sqrt{n}\g\|\bx_k\!-\!\bx_{k-\!1}\|,
	\end{align*}
where we use the gradient Lipschitz continuity property for each $f_i$, the $s$-update in \eqref{eq-s} the definitions of the $D$-quantity in~\eqref{eq-x-D-quants} and the constant $\varphi_k$ in~\eqref{eq-const-all}. Hence, we obtain the following relation for the second term in \eqref{eq-proofweightx}:
	\begin{align*}
    &\sum_{i=1}^n[\phi_{k+1}]_i\|y_k^i-n [\pi_k]_i \nabla f(\hat x_k)\|
    \\
    \le& S(\by_k,\pi_k) +L\sqrt{n}\varphi_kD(\bx_k,\phi_k)
    +L\sqrt{n}\g\|\bx_k-\bx_{k-1}\|.
    \end{align*}
	For the final term in~\eqref{eq-proofweightx}, since $\phi_k$ is stochastic, we have 
	\[\sum_{i=1}^n[\phi_{k}]_i\|x_k^i-x_{k-\!1}^i\| \!\le\! \sqrt{\!\sum_{i=1}^n[\phi_{k}]_i\|x_k^i-\!x^i_{k\!-\!1}\|^2} \!\le \! \|\bx_k\!-\!\bx_{k-\!1}\|.\]
	We can obtain similar relation when the weight is $\phi_{k+1}$, thus,
	\[\sum_{i=1}^n (\b[\phi_{k+1}]_i+\g[\phi_{k}]_i)\|x_k^i-x_{k-1}^i\| \le (\b+\g)\|\bx_k-\bx_{k-1}\|.\]
	Substituting the estimates obtained above for each term on the right-hand-side of~\eqref{eq-proofweightx} gives the desired relation.
\end{proof}

\subsection{Proof of Proposition~\ref{prop-xcontract}}\label{app:Prop2}
\begin{proof}
	Let $u_k^i=x_k^i-x_{k-1}^i$ and $\hat u_k=\sum_{i=1}^n[\phi_k]_i u_k^i$. Subtracting the relation for $\hat x_k$ given in Lemma~\ref{lem-weightedavg} from the $x$-update in equation~\eqref{eq-xcomb}, and using the triangle inequality, we have 
	\begin{align}\label{eq-11}
	D(\bx_{k+1}\!,\phi_{k+1})
	\le \!&\sqrt{\sum_{i=1}^n [\phi_{k+1}]_i \Bigg\|\sum_{j=1}^n [A_k]_{ij} x_k^j -\hat x_k\Bigg\|^2} \nonumber\\
	+\a 
	&\sqrt{\sum_{i=1}^n [\phi_{k+1}]_i \Bigg\|y_k^i \!- \!\!\sum_{j=1}^n[\phi_{k+1}]_j y_k^j\Bigg\|^2}\nonumber\\
	+\g&\sqrt{\sum_{i=1}^n [\phi_{k+1}]_i \Bigg\|\sum_{j=1}^n [A_k]_{ij} u_k^j -\hat u_k\Bigg\|^2}\nonumber\\
	+\b &\sqrt{\sum_{i=1}^n [\phi_{k+1}]_i \Bigg\|u_k^i \!- \!\!\sum_{j=1}^n[\phi_{k+1}]_j u_k^j\Bigg\|^2}.
	\end{align}
	Under Assumption~\ref{asm-graphs}, we use Lemma~\ref{lem-basic-xcontract} to estimate the first term in \eqref{eq-11}, with $A=A_k$, $x_i=x_k^i$ and 
	$\phi_{k+1}^{\T}A_k=\phi_k^{\T}$:
	\begin{align}\label{eq-proof-1st}
	\sqrt{\!\sum_{i=1}^n [\phi_{k+1}]_i \Bigg\|\sum_{j=1}^n [A_k]_{ij} x_k^j \!-\!\hat x_k\Bigg\|^2}
	\!&\le\! c_k
	\sqrt{\!\sum_{j=1}^n [\phi_k]_j \|x_k^j \!-\! \hat x_k\|^2} \nonumber\\
	&\le c_k
	D(\bx_k,\phi_k).
	\end{align}
	Similarly for the third term in \eqref{eq-11}, using Lemma~\ref{lem-basic-xcontract} with  $A=A_k$, $x_i=u_k^i$ and 
	$\phi_{k+1}^{\T}A_k=\phi_k^{\T}$, we obtain
	\begin{align*}
	\sqrt{\!\sum_{i=1}^n [\phi_{k+1}]_i \Bigg\|\sum_{j=1}^n [A_k]_{ij} u_k^j \!-\!\hat u_k\Bigg\|^2}
	\!\le\! c_k
	\sqrt{\!\sum_{j=1}^n [\phi_k]_j \|u_k^j \!-\! \hat u_k\|^2}.
	\end{align*}
    Next, using the relation in Lemma~\ref{lem-normlincomb} with $\g_i=[\phi_{k}]_i$, $u_i=u_k^i$ and $u=0$,
	it follows that
	\[\sum_{i=1}^n [\phi_{k}]_i \Bigg\|u_k^i - \sum_{j=1}^n[\phi_{k}]_j u_k^j\Bigg\|^2
	\!\!\le\!\! \sum_{i=1}^n [\phi_{k}]_i \|u_k^i \|^2 \le \|\bx_k-\bx_{k-1} \|^2\!\!,\]
	where the last inequality follows from the stochasticity of $\phi_k$ and the definition of $u_k^i$, for all $i\in [n]$. Thus,
	\begin{align} \label{eq-proof-3rd}
	\sqrt{\!\sum_{i=1}^n [\phi_{k+1}]_i \Bigg\|\sum_{j=1}^n [A_k]_{ij} u_k^j \!-\!\hat u_k\Bigg\|^2}
	\!\le\! c_k
	\|\bx_k-\bx_{k-1} \|.
	\end{align}
	Similar argument can be used to estimate the forth term in \eqref{eq-11}, which yields
	\begin{align} \label{eq-proof-4th}
	&\sum_{i=1}^n [\phi_{k+1}]_i \Bigg\|u_k^i - \sum_{j=1}^n[\phi_{k+1}]_j u_k^j\Bigg\|^2 \nonumber \\
	\le & \sum_{i=1}^n [\phi_{k+1}]_i \|u_k^i \|^2 \le \|\bx_k-\bx_{k-1} \|^2,
	\end{align}
	and, for the second term in \eqref{eq-11}, as follows
    \begin{align}\label{eq-proof-2nd}
    &\sqrt{\sum_{i=1}^n [\phi_{k+1}]_i \Bigg\|y_k^i - \sum_{j=1}^n[\phi_{k+1}]_j y_k^j\Bigg\|^2}
	\!\le\! \sqrt{\sum_{i=1}^n [\phi_{k+1}]_i \|y_k^i \|^2}\nonumber\\
	\le&\sqrt{\max_{j\in[n]} ([\phi_{k+1}]_j[\pi_k]_j)} \sqrt{\sum_{i=1}^n\frac{\left\|y_k^i \right\|^2}{[\pi_k]_i}}
	\le \sqrt{\sum_{i=1}^n\frac{\left\|y_k^i \right\|^2}{[\pi_k]_i}} \nonumber\\
	\le &L\sqrt{n} \left(\varphi_k\|\hat x_k \!-\! x^*\| \!+\! \varphi_kD(\bx_k,\!\phi_k)\!+\!\g\|\bx_k \!-\!\bx_{k-1}\|\right)\nonumber\\
    &\!+\! S(\by_k,\!\pi_k),
	\end{align}
	where the last inequality follows from Lemmas~\ref{lem-yinv} and~\ref{lem-ysum}.
By combining the estimates for each term in \eqref{eq-11}, as given by relations \eqref{eq-proof-1st}--\eqref{eq-proof-2nd}, we arrive at the desired relation.
\end{proof}

\subsection{Proof of Proposition~\ref{prop-ycontract}}\label{app:Prop4}
\begin{proof}
	By defining $w_k^i=\sum_{j=1}^n [B_k]_{ij} y_k^j$, $\bw_k=(w_k^1,\ldots,w_k^n)$ and $\bg_k=(\nabla f_1(s_k^1),\ldots,\nabla f_n(s_k^n))$, the update for the $y$-iterate in compact form is given as
	\begin{equation}\label{eq-ycomp}
	\by_{k+1}=\bw_k+ \bg_{k+1}-\bg_k\qquad\hbox{for all }k\ge0.
	\end{equation}
	Let $\Lambda=\diag^{-1}(\pi_{k+1})$, we can write
	\begin{equation*}
	\by_{k+1}\Lambda=\bw_k\Lambda
	+ (\bg_{k+1}-\bg_k)\Lambda\qquad\hbox{for all }k\ge0.
	\end{equation*}
	By subtracting the vector 
	$\bar{\by}_{k+1} =(\bar{y}_{k+1},\ldots,\bar{y}_{k+1})$, where $\bar{y}_{k+1}=\sum_{j=1}^n y_{k+1}^j$,
	from both sides of the preceding relation, we have for all $k\ge0,$
	\begin{equation*}
	\by_{k+1}\Lambda-\bar{\by}_{k+1}  
	\!=\bw_k\Lambda-\bar{\by}_k +(\bar{\by}_{k} \!-\bar{\by}_{k+1} )+ (\bg_{k+1}\!-\bg_k)\Lambda.
	\end{equation*}
	By taking $\pi_{k+1}$-induced norm on both sides and noting that $S(\by_{k+1},\pi_{k+1})=|\by_{k+1}\Lambda-\bar{\by}_{k+1} |{\pi{k+1}}$, we obtain
	\begin{align}\label{eq-srel1}
	S(\by_{k+1},\pi_{k+1})&\le \|\bw_k\Lambda-\bar{\by}_k\|_{\pi_{k+1}}+\|\bar{\by}_{k+1} -\bar{\by}_{k} \|_{\pi_{k+1}}\!\nonumber\\
	&+\|(\bg_{k+1}-\bg_k)\Lambda\|_{\pi_{k+1}}.
	\end{align}
	For the first term in \eqref{eq-srel1}, by using the definitions of $\bw_k$ and $\bar{\by}_k$, we can deduce that 
	\begin{align*}
	\|\bw_k\Lambda-\bar{\by}_k\|_{\pi_{k+1}}
	=\sqrt{\sum_{i=1}^n[\pi_{k+1}]_i \left\|\frac{w_k^i}{[\pi_{k+1}]_i}-\sum_{\ell=1}^n y_k^\ell\right\|^2}\\
	\le \tau_k 
	\sqrt{\sum_{i=1}^n [\pi_k]_i \left\|\frac{y_i}{[\pi_k]_i} - \sum_{\ell=1}^n y_\ell\right\|^2} = \tau_k\,S(\by_k,\pi_k),
	\end{align*}
	where the first inequality follows from Lemma~\ref{lem-basic-ycontract}
	with $\bbG=\bbG_k$, a strongly connected directed graph (see Assumption~\ref{asm-graphs}), $B=B_k$, $\pi=\pi_{k+1}$, and $\nu=\pi_k$. The last equality follows from the definition of $S(\by_k,\pi_k)$.
 
	For the second term in \eqref{eq-srel1}, since $\bar{\by}_{k} =\sum_{i=1}^n y_k^i=\sum_{i=1}^n \nabla f_i(s_k^i)$ (as stated in Lemma~\ref{lem-sumgrad}), we have
	\begin{align*} 
	&\|\bar{\by}_{k+1} \!-\!\bar{\by}_{k} \|_{\pi_{k+1}}\!=\!\sqrt{\!\sum_{i=1}^n [\pi_{k+1}]_i\|\bar{y}_{k+1}\!-\!\bar{y}_k\|^2}
	\!=\!\|\bar{y}_{k+1}\!-\!\bar{y}_k\|\\
	&=\!\Bigg\|\!\sum_{i=1}^n\! \left(\nabla \!f_i(s_{k+1}^i)\!-\! \nabla\! f_i(s_k^i) \right)\!\!\Bigg\|
	\!\le\! \sum_{i=1}^n \!\|\nabla\! f_i(s_{k+1}^i)\!-\! \nabla\! f_i(s_k^i)\|\\
	&\le L\sum_{i=1}^n \|s_{k+1}^i - s_k^i\|\le L\sqrt{n}\|\bs_{k+1}-\bs_k\|,
	\end{align*}
where the last inequality follows from
the Lipschitz continuity of the gradients $\nabla f_i$ (Assumption~\ref{asm-functions}).

	For the last term in relation~\eqref{eq-srel1}, we have
	\begin{align*}\|(\bg_{k+1}-\bg_k)\Lambda\|_{\pi_{k+1}}
	=\sqrt{\sum_{i=1}^n   \frac{\|\nabla f_i(s_{k+1}^i)- \nabla f_i(s_k^i)\|^2}{[\pi_{k+1}]_i}} \\
	\le L\sqrt{\sum_{i=1}^n \frac{\|s_{k+1}^i - s_k^i\|^2}{[\pi_{k+1}]_i} }
	\le \frac{L}{\sqrt{\min(\pi_{k+1})}}\,\|\bs_{k+1}-\bs_k\|,
	\end{align*}
	where the first inequality follows by the Lipschitz continuity of the gradients $\nabla f_i$. Thus,
	\begin{align*}
	S(\by_{k+1},\pi_{k+1})
	\le \tau_k\, S(\by_k,\pi_k)+L r_k\|\bs_{k+1}-\bs_k\|,
	\end{align*}
	where $r_k=\sqrt{n}+ \frac{1}{\sqrt{\min(\pi_{k+1})}}$. Using the compact form of the $s$-update of the method in~\eqref{eq-s}, we further obtain
	\begin{align*}
	\|\bs_{k+1}-\bs_k\| \le (1+\g)\|\bx_{k+1}-\bx_k\|+\g\|\bx_{k}-\bx_{k-1}\|.
	\end{align*}
	The desired relation follows from the preceding two relations.
\end{proof}

\vfill
	
\end{document}